\documentclass[11pt,a4paper,oneside]{article}
\usepackage[metapost]{mfpic}
\usepackage{subfigure}
\usepackage{graphicx}
\usepackage{algorithm,algorithmic}
\usepackage[dvipsnames]{color}
\usepackage{amsmath}
\usepackage{amsthm}

\usepackage[bookmarks=false]{hyperref}

\usepackage{geometry}
\usepackage{color}
\usepackage{url}
\usepackage[font=small,labelfont=bf, figurename=Fig.]{caption}

\usepackage{fancyhdr} 

\usepackage{graphicx}

\definecolor{hellgrau}{gray}{.8}
\definecolor{dunkelblau}{rgb}{0, 0, .7}
\definecolor{roetlich}{rgb}{1, .7, .7}
\definecolor{dunkelmagenta}{rgb}{.3, 0, .3}

\def \4th{{1 \over 4}}

\def \2th{{1 \over 2}}

\normalsize
\title{Leveraging NMF to Investigate Air Quality in Central Taiwan}

\begin{document}
\pagestyle{fancy}



\author{ Shu-Chuan~Chen$^{a}$, Jui-Fang~Chang$^{b}$,  ~and ~
Yintzer Shih$^{b}$\footnote{Corresponding author. Supported by the National Science Technology
Council of Taiwan through Project NSTC 112-2115-M-005-003.}\\[4mm]
$^a$ Department of Mathematics and Statistics, Idaho State University, Pocatello, Idaho, USA\\[1mm]
$^b$ Department of Applied Mathematics, National Chung Hsing University, 
Taichung, Taiwan\\[1mm]
}

\date{}
\maketitle
\vspace{1mm}

\begin{abstract}
This study investigates air pollution in central Taiwan, focusing on key pollutants, including SO$_2$, NO$_2$, PM$_{10}$, and PM$_{2.5}$. We use non-negative matrix factorization (NMF) to reduce data dimensionality, followed by wind direction analysis and speed to trace pollution sources. Our findings indicate that PM$_{2.5}$ and NO$_2$ levels are primarily influenced by local sources, while SO$_2$ levels are more affected by transboundary factors. For PM$_{10}$, contributions from domestic and transboundary sources are nearly equal.\end{abstract}


\noindent


\section{Introduction}

Throughout history, the rapid advancement of human civilization and technology has led to significant environmental and air pollution challenges. The widespread use of transportation fuels and various energy sources has generated large quantities of pollutants. Major pollution events—such as the 1940 Los Angeles smog incident, the 1952 London Great Smog, the 1961 Yokkaichi asthma incident in Japan, and the 2013 air quality crisis in Mainland China—have been directly linked to high energy consumption. These events have posed serious risks to human health, leading to various diseases, environmental damage, and a decline in the quality of life.

The Taichung Thermal Power Plant (TTPP) in Taichung City has consistently drawn public and political attention, particularly during election cycles from 2009 to 2019. As Taipower’s largest thermal power plant, TTPP has a total installed capacity of 5.788 million kilowatts, representing approximately 20\% of Taiwan’s total power generation. From 2006 to February 2017, it was the world’s largest thermal power plant. This prominence highlights the importance of studying air pollution in central Taiwan, especially within the Greater Taichung City area.

Public concern about air pollution is growing due to its harmful effects on health and the environment. Polluted air can damage the respiratory, cardiovascular, and nervous systems, with air pollutants taking various forms--solid particles, liquid droplets, and gases--each impacting health differently. A significant portion of atmospheric sulfur dioxide (SO$_2$) and nitrogen dioxide (NO$_2$) arises from fuel combustion in power plants, vehicles, and industries, although the contributions from each source vary. Both SO2 and NO2 can irritate the eyes, mouth, nose, and respiratory tract, and they are major contributors to acid rain. 

Particulate matter (PM) is categorized by particle size, primarily into PM$_{10}$, which includes particles 10 micrometers or smaller, and PM$_{2.5}$, which consists of finer particles 2.5 micrometers or smaller. SO$_2$ generally remains in the troposphere for about two days, or roughly one day in areas with high moisture at cloud tops. In central Taiwan, the mountain range disrupts prevailing winds, causing air pollution in inland cities and the southwest region largely from local sources. Kishcha et al. highlight that air pollution is directly associated with aerosol optical depth (AOD), which affects atmospheric water circulation, cloud formation, precipitation, and climate, demonstrating the broad indirect effects of pollution.

In recent years, PM$_{2.5}$ has become a research focus due to its small size--smaller than bacteria--making it easily inhalable and potentially harmful to health. In 2015, PM$_{2.5}$ was identified as the leading contributor to the Global Burden of Disease 
\cite{long}. To address air pollution effectively, it is essential first to identify primary pollution sources and then develop practical, targeted solutions.

In Taiwan, air pollution sources are categorized as either transboundary or domestic. Lin et al. \cite{wind} note that transboundary pollution mainly originates from the Asian continent, especially in winter when the northeast monsoon prevails. Key sources include dust storms, frontal pollution, and background air masses. According to Li et al. \cite{1-2}, PM$_{2.5}$ from transboundary pollution is largely sourced from North China and the East China Sea, with dust storms and northeast monsoons in spring and winter transporting Asian dust and artificial particles across the Taiwan Strait. This paper examines the main sources of pollutants such as O$_3$, NO$_2$, PM$_{10}$, and PM$_{2.5}$ in central Taiwan.

The paper is structured as follows: Section 2 outlines the methods, Section 3 describes data sources, Section 4 presents the analysis results, followed by a validation of the proposed method in Section 5. Section 6 concludes the study.

\section{Methods and Algorithms}
\label{set2}

 
Non-negative matrix factorization (NMF), introduced by Lee and Seung in 1999 
\cite{nmf3}, has become a valuable tool across various fields. A key advantage of NMF is that all values in the resulting matrices are non-negative, which is especially beneficial in practical applications. In scientific research, matrix factorization is commonly used for efficient data processing, and NMF offers a simpler approach with easily interpretable results compared to traditional methods. Techniques such as PCA, SVD (Singular Value Decomposition), and QR factorization are often used for matrix factorization; however, they can produce negative values in the output even if all input elements are positive. Although these negative values may be mathematically valid, they are often irrelevant or impractical in real-world scenarios. For example, in image processing, where pixel values range from 0 to 255, negative values in the decomposition lack practical meaning. NMF effectively circumvents this issue by ensuring that all outputs are non-negative. Additionally, NMF has broad applications, including in speech processing, robotics, and fields such as biomedical and chemical engineering.

\subsection{Nonnegative Matrix Factorization}

We now introduce the NMF algorithm, following the approach outlined by Lee \cite{nmf1,nmf2}. First, we represent the data as a matrix  $A$ and apply the NMF method to decompose  $A$ into two matrices, $W$ and $H$, as shown in the following equation:
 \begin{equation} A_{m \times n} \simeq W_{m \times k} \times H_{k \times n} = \hat{A}. \label{nmf} \end{equation}
The non-negative initial matrices $W$ and $H$ are randomly initialized, and we use Euclidean distance to calculate the cost function, as shown in the following formula:
\begin{equation}
J(W,H)=\frac{1}{2}\Vert A-WH\Vert^2_{F}=\frac{1}{2}\sum_{ij}(A_{ij}-(WH)_{ij})^2.
\label{cost}
\end{equation}
The multiplicative update rules are implemented to iterate until convergence, as shown in the following formula:
\begin{equation}
W_{il}=W_{il}\cdot \frac{(AH^T)_{il}}{(WHH^T)_{il}} , ~~H_{lj}=H_{lj}\cdot \frac{(W^TA)_{lj}}{(W^TWH)_{lj}}.
\label{rules}
\end{equation}

The formula (\ref{rules}) ensures that the matrix elements are all positive values. We will now prove the multiplication update algorithm. We derive it from matrix multiplication and differentiation:
\begin{equation}
(WH)_{ij}=\sum_{l}W_{il}H_{lj} ~~\mbox{and}~~ \frac{\partial(WH)_{ij}}{\partial W_{il}}=H_{lj}.
\label{update}
\end{equation}
Then
\begin{equation}
\begin{split}
\frac{\partial }{\partial W} J(W,H)\bigg|_{il} &= \sum_{j}[-H_{lj}(A_{ij}-(WH)_{ij})]H_{lj} \\ &= \sum_{j} -A_{ij}H_{lj}+\sum_{j}(WH)_{ij}H_{lj} \\ &= -(AH^T)_{il}+(WHH^T)_{il}.
\end{split}
\label{partial}
\end{equation}
Following \eqref{partial}, we have
\begin{equation}
\frac{\partial }{\partial H}J(W,H)\bigg|_{lj}=-(W^TA)_{lj}+(W^TWH)_{lj}.
\label{partial1}
\end{equation}
Next, using the gradient, we have
\[ W_{il}\Leftarrow W_{il}-\alpha_{1} \cdot [-(AH^T)_{il}+(WHH^T)_{il}]\text{,} \]
\[ H_{lj}\Leftarrow H_{lj}-\alpha_{2} \cdot [-(W^TA)_{lj}+(W^TWH)_{lj}]. \]
Let
\[ \alpha_{1}=\frac{W_{il}}{(WHH^T)_{il}}\text{,}  \quad \alpha_{2}=\frac{H_{lj}}{(W^TWH)_{lj}}.\]
Then, we have the result of equation \eqref{rules}.
Finally, one normalizes by converting the data to values between $0$ and $1$, and the formula is as follows:
\begin{equation}
H_{ij}=
\begin{cases}
\frac{H_{ij}-H_{\min}}{H_{\max}-H_{min}}, & \mbox{if $H_{\min} < H_{\max}$}\\
1, &  \mbox{if $H_{\min} = H_{\max}$}.
\end{cases}
\label{nor}
\end{equation}

The NMF algorithm is described below.

\begin{algorithm}
\caption{NMF}
\label{alg:NMF}
\begin{algorithmic}
\STATE {Input:
$A, k $, where $A$ is the data matrix, $k$ is the required dimension reduction.}
\STATE {Step 1. Given an initial nonnegative matrices $W$ and $H$,  normalize every row in matrix $H$.}
\STATE {Step 2. While (not converged) \\
                   \hspace{1.2cm} for $i=1:m$, ~~~for $l=1:k$
               \[ W_{il} \leftarrow W_{il}\cdot \frac{(AH^T)_{il}}{(WHH^T)_{il}} \]                         
                  \hspace{1.2cm} for $j=1:n$,  ~~~for $l=1:k$\\
               \[  H_{lj}  \leftarrow H_{lj}\cdot \frac{(W^TA)_{lj}}{(W^TWH)_{lj}} \]  
                \hspace{1.2cm} Renormalize each row in $H$}\\
                \hspace{1.2cm} End
\STATE {Output: Matrices $W$ and $H$.}
\end{algorithmic}
\end{algorithm}

\section{Data sources}

The Environmental Protection Department provided the research data over the years \cite{EPAA}, consisting of hourly pollution concentration values. Our study focuses on data collected from air quality monitoring stations across 14 regions in Taiwan, including Miaoli, Taichung, Changhua, and Nantou. The main pollutants analyzed are SO$_2$ (sulfur dioxide), NO$_2$ (nitrogen dioxide), PM$_{2.5}$, and PM$_{10}$. These pollutants are measured in parts per billion (ppb) for SO$_2$ and NO$_2$, and in $\mu g/m^3$ for PM$_{2.5}$ and PM$_{10}$. The study period spans from 2008 to 2017, encompassing 87,672 records. The data are stored in an 87,672 $\times$ 14 matrix for each pollutant, resulting in four matrices of the same dimensions. 

Fig. \ref{M} illustrates the locations of the 14 monitoring stations. Industrial air quality monitoring stations are positioned downwind of industrial areas to assess the impact of industrial pollution. General air quality monitoring stations are situated in densely populated regions where high pollution levels may occur, while the background air quality monitoring station is placed in areas with minimal synthetic pollution \cite{EPA}.

\bgroup-
\begin{figure}[hp]
	\begin{center}
		\begin{center}
			\includegraphics[height=8cm]{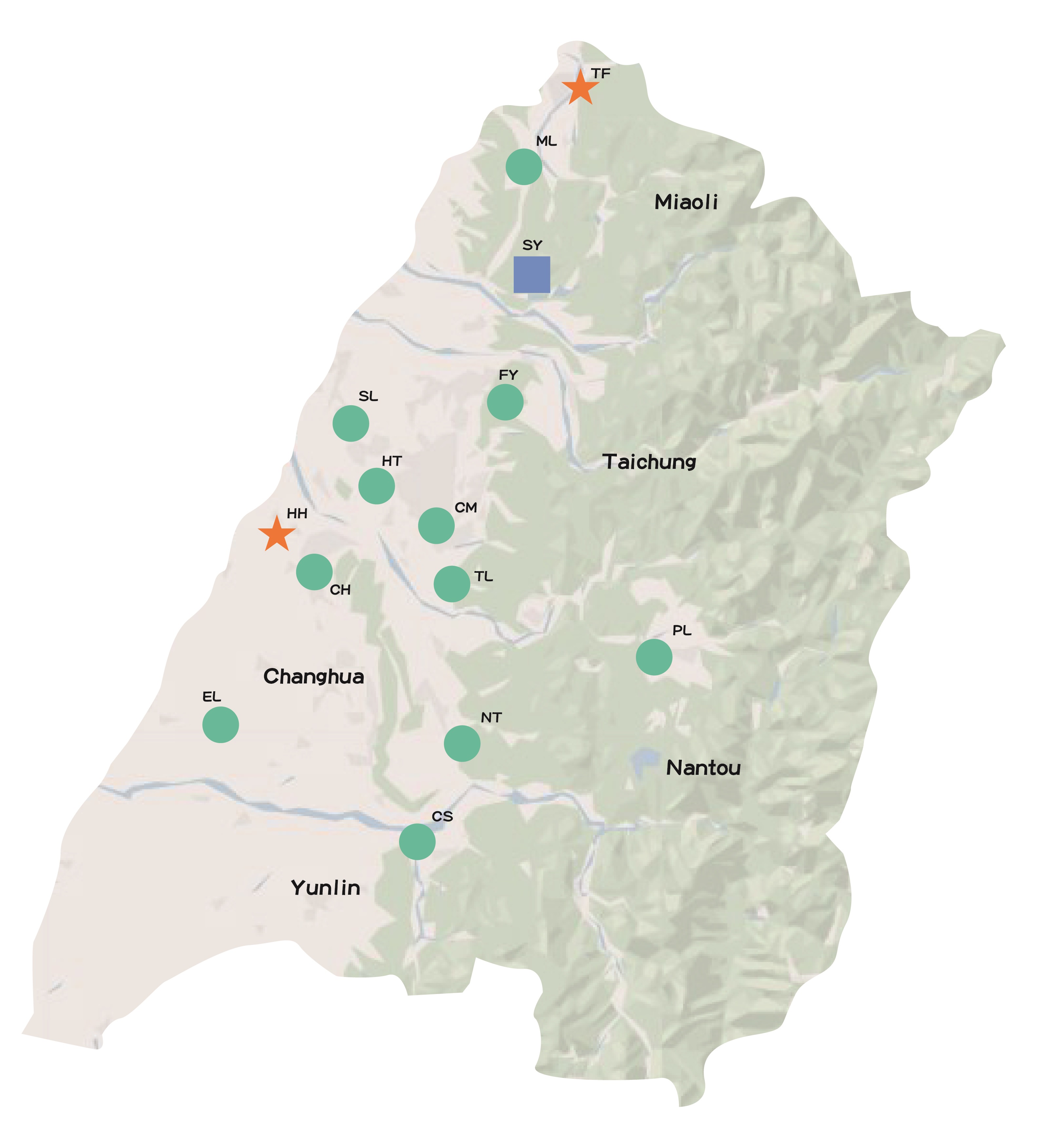}
		\end{center}
		\caption{Central Taiwan Terrain Map: Stars indicate Industrial Air Quality Monitoring Stations; squares indicate Background Air Quality Monitoring Stations; and circles indicate General Air Quality Monitoring Stations.}
		\label{M}
	\end{center}
\end{figure}
\egroup

We visualized the concentrations of SO$_{2}$, NO$_{2}$, PM$_{10}$, and PM$_{2.5}$, and divided the average concentrations of air pollution data into hourly and monthly intervals (Fig. \ref{2.1}-\ref{2.2}). First, the monitoring stations with the highest concentrations of NO$_{2}$ (Fig. \ref{2.1}) are CM, TL, CH, and NT, with CM being the most polluted. Higher NO$_{2}$ concentrations are observed between 5:00 PM and 12:00 AM at all monitoring stations. Second, the stations with the highest concentrations of SO$_{2}$ are HH and CH in Changhua, with HH having the highest levels. The HH station is an industrial air quality monitoring station positioned downwind of the thermal power plant in Taichung and influenced by the prevailing northeast wind. Third, PM$_{10}$ concentrations rise after noon, especially at the EL, NT, CS, HT, and TL stations. Fourth, the stations with the lowest PM$_{2.5}$ concentrations are TF, ML, and SY, all located in Miaoli. The particulate matter levels in Miaoli are similar to those in other regions, but overall, the SY station has the best air quality in central Taiwan.

In Fig. \ref{2.2}, we observe that the highest pollutant concentrations occur in March, while the lowest concentrations of NO$_{2}$, PM$_{10}$, and PM$_{2.5}$ are seen during the summer. SO$_{2}$ concentration remains relatively consistent throughout the year. Additionally, seasonal data (Table \ref{bmonth}) show that pollution is more severe during the spring and winter.
Wind speed, a key factor in our analysis, will be further explored later. Following \cite{wind}, wind speed is classified as shown in Table \ref{speed}.

\bgroup
\begin{figure}[ph]
	\begin{center}
		\begin{center}
			\includegraphics[width=12cm]{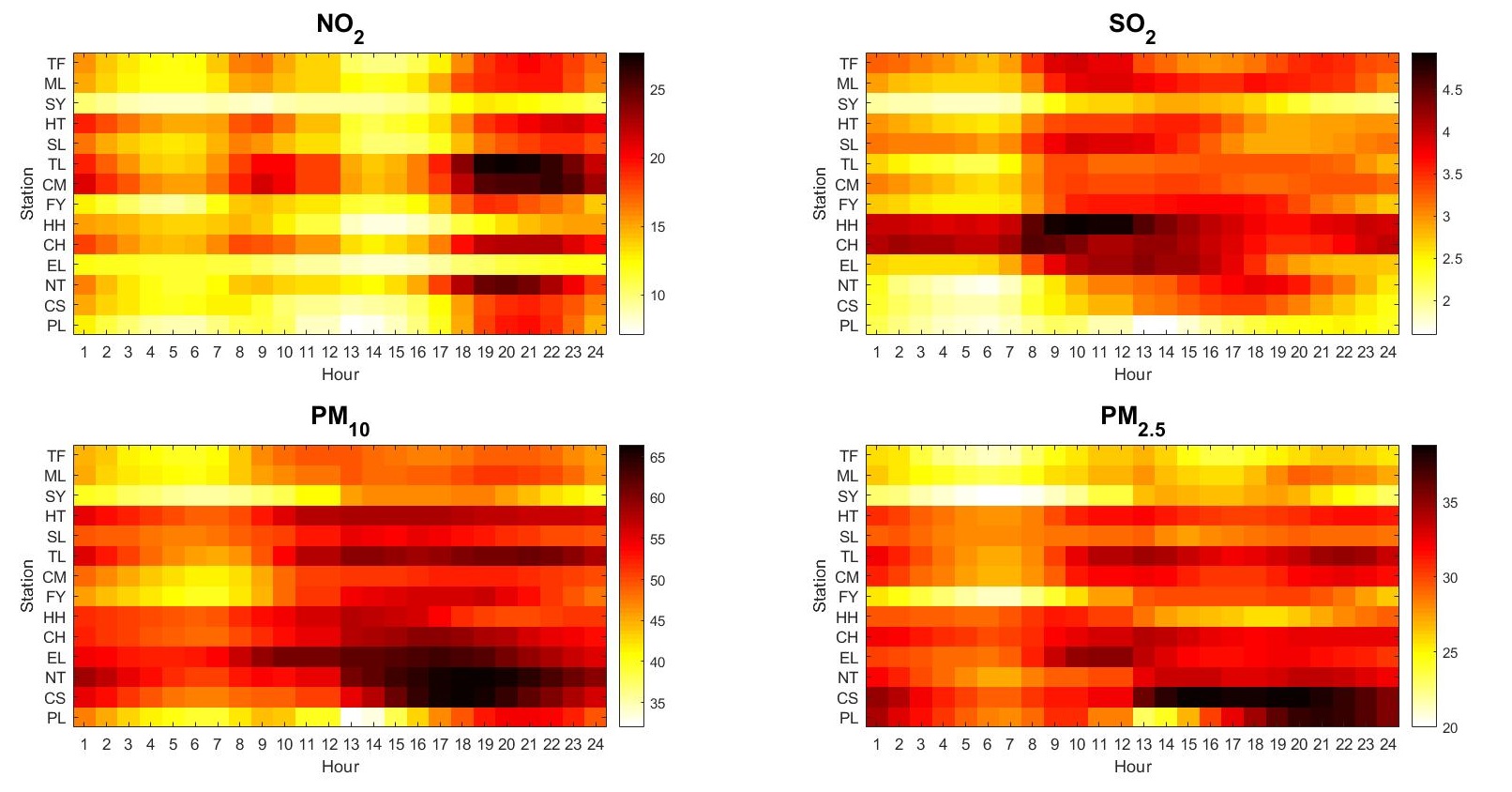}
		\end{center}
		\caption{Visualization for hourly average concentration}
		\label{2.1}
	\end{center}
\end{figure}
\egroup

\bgroup
\begin{figure}[ph!]
	\begin{center}
		\begin{center}
			\includegraphics[width=12cm]{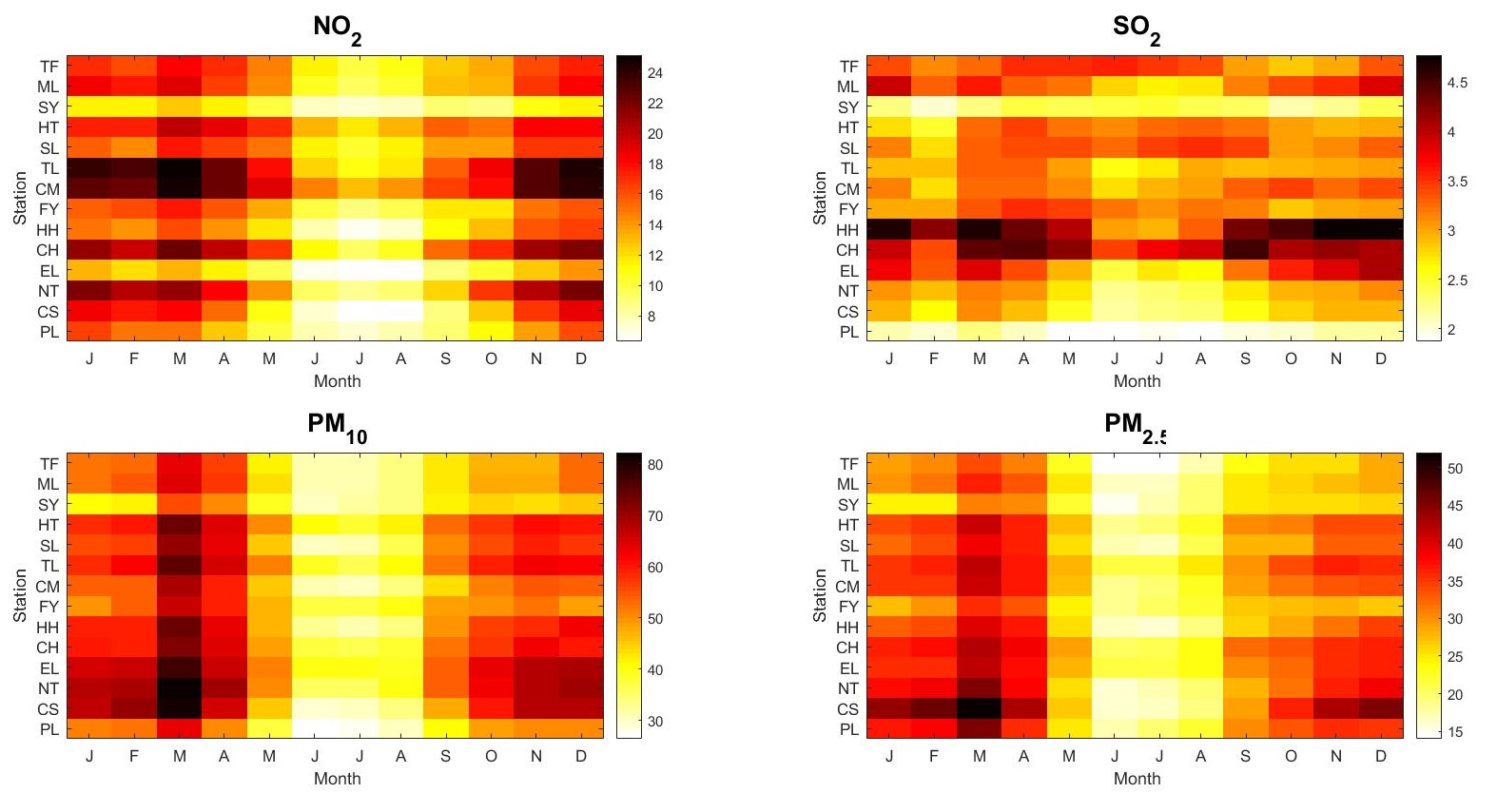}
		\end{center}
		\caption{Visualization for monthly average concentration.}
		\label{2.2}
	\end{center}
\end{figure}
\egroup

\bgroup
\begin{table}[ph]
	\begin{center}
		\caption{Wind speed}
		\label{speed}
		\begin{tabular}{|c|c|}
		\hline
speed ($m/s$)  & definition \\ \hline
$<0.2$  & Calm  \\ \hline
$0.2  \sim 5.4$ & Gentle breeze  \\ \hline
$5.5 \sim 13.8$  & Strong   \\ \hline
$\ge 13.9$  & Gale   \\ \hline
		\end{tabular}
	\end{center}
\end{table}
\egroup

\bgroup
\begin{table}[h]
	\begin{center}
		\caption{The percentage of the four-season concentration}
		\label{bmonth}
		\begin{tabular}{|c|c|c|c|c|}
		\hline
proportions & spring & summer  & autumn  & winter \\ 
 &  (Mar.-May) &  (Jun.-Aug.) &  (Sep.-Nov.) &  (Dec.-Feb.) \\ \hline
 NO$_{2}$  & 28.23\% & 16.6\%	& 24.83\% & 30.33\% \\ \hline
 SO$_{2}$  & 26.27\% &	22.98\%	& 25.37\%	& 25.37\% \\ \hline
 PM$_{2.5}$  & 28.99\%	& 16.24\%	&25.81\%	& 28.96\%  \\ \hline
 PM$_{10}$  & 28.95\% &	16.92\%	& 26.03\%	& 28.1\%  \\ \hline
		\end{tabular}
	\end{center}
\end{table}
\egroup

\section{NMF analysis}

Fig. \ref{F} illustrates a flowchart detailing the process. The initial step involves identifying the optimal value of $k$, which yields $k$ features. Subsequently, by utilizing data on wind direction, wind speed, and pollution concentration, we ascertain whether each feature is linked to domestic or foreign sources of pollution.

\subsection{Choose $k$}

Prior to conducting the analysis, it is imperative to ascertain the appropriate value of k. For this purpose, we employ the method proposed by Jean et al. \cite{jean}, wherein we select the value of k at which the Cophenetic Correlation begins to exhibit a decline. According to Limem et al. \cite{Limem}, the value of k must adhere to the stipulation $nk + km \ll nm$, which simplifies to $k \ll \min\{ n, m \}$. Given that $k \ll 14$, we consider values within the range of $2 \leq k \leq 7$. 

In Figs. \ref{N}(a) and \ref{N}(b), we analyze the scenarios for $k=6$ and $k=5$, respectively. For the case where $k=5$, the correlation exhibits a significant decline to a minimum point, followed by an increase, thereby suggesting an optimal clustering effect. In Figs. \ref{N}(c) and \ref{N}(d), we observe the cases for $k=4$ and $k=3$, respectively. Both cases reveal multiple points of decline; hence, we opt for the scenario that demonstrates the most pronounced decrease to achieve optimal clustering.

\bgroup
 \begin{figure}[h]
	\begin{center}
		\begin{center}
			\includegraphics[height=10cm]{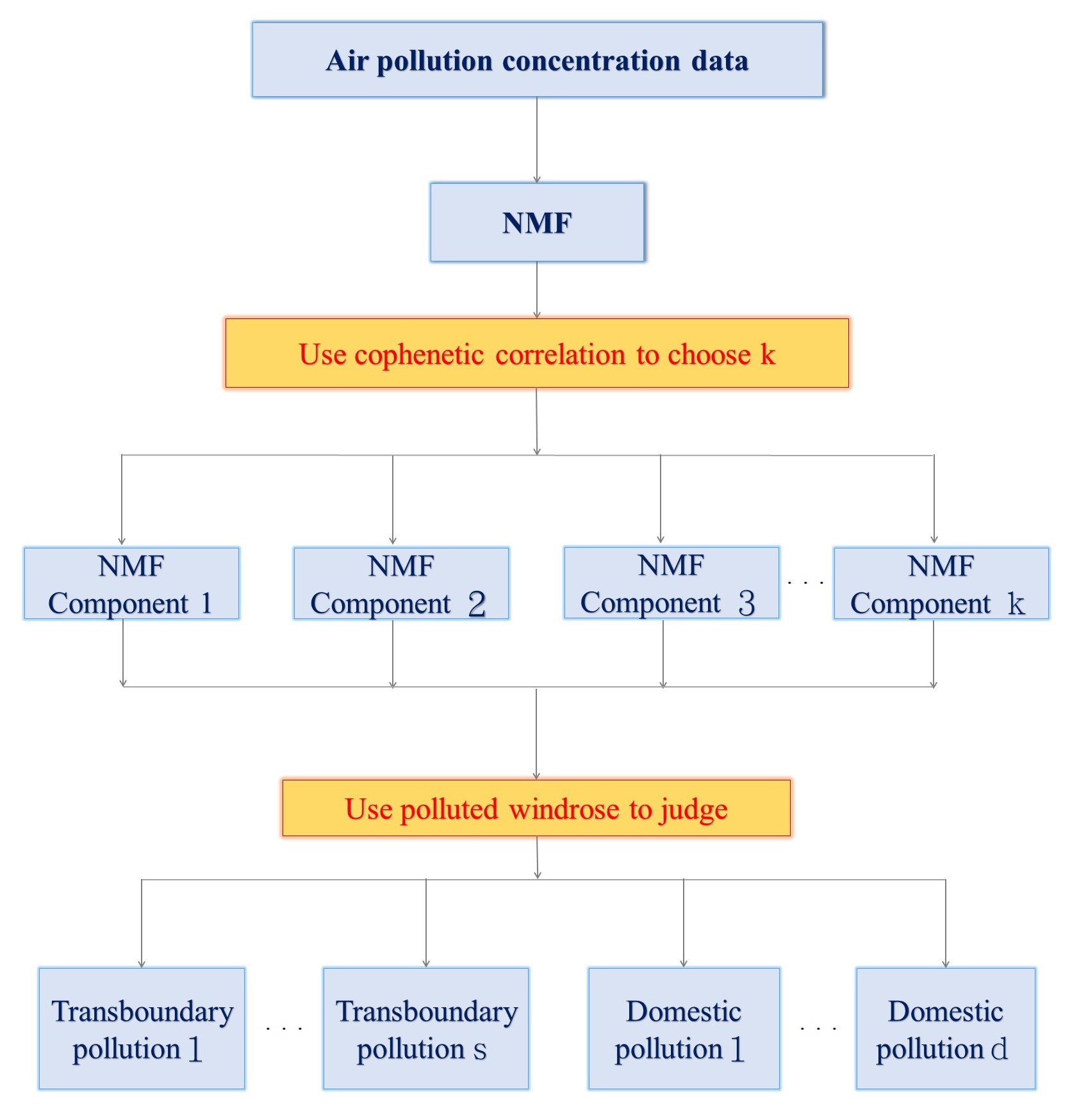}
		\end{center}
		\caption{The flow chart for proposed data analysis}
		\label{F}
	\end{center}
\end{figure}
\egroup

\bgroup
 \begin{figure}[h]
	\begin{center}
		\begin{center}
			\includegraphics[height=8cm]{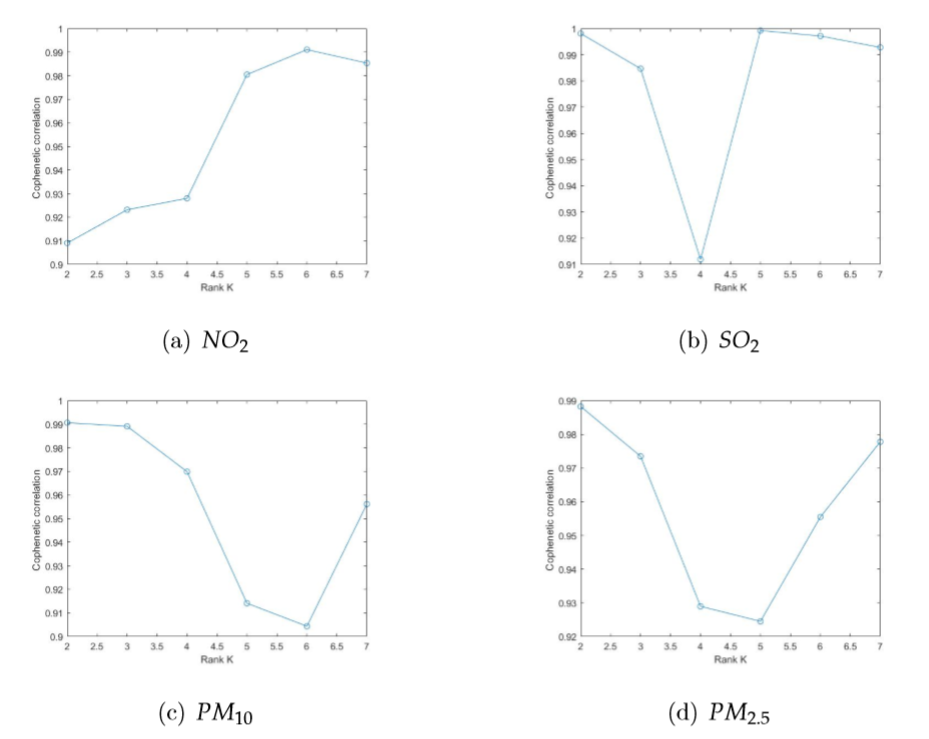}
		\end{center}
		\caption{NMF for cophenetic correlation}
		\label{N}
	\end{center}
\end{figure}
\egroup

\section{Results}

We refer to the work by Luoa et al. \cite{NMF}. First, matrix A is decomposed using the NMF method. Following the principles outlined in Section \ref{set2}, the resulting matrices W and H are then decomposed and visualized. 
 In Fig. \ref{hno2}, we observe that the monitoring stations with the highest contributions to NO$_{2}$ concentration in NMF1 are NT, CS, and PL, all located in Nantou. Fig. \ref{wno2} illustrates significant seasonal variations, with severe pollution levels occurring in the evening at NMF1. Additionally, NMF1 is notably influenced by westerly winds, with high pollution levels concentrated at a wind speed of 4 m/s, as shown in Fig. \ref{RNO2}. Based on this, we classify NMF1 as domestic pollution, contributing 16.79\% of all features.

NMF2 significantly contributes to pollutant concentrations at stations such as TL and CH, peaking during periods of heavy traffic flow. Although strong northerly winds are present during autumn and winter, the concentrations remain low. However, high pollution levels are concentrated at a wind speed of 3.5 $m/s$ throughout all seasons. Thus, we classify NMF2 as domestic pollution, contributing 16.58\% of all features.

For NMF3, the monitoring stations with high contributions to pollutant concentration are HH and EL in Changhua. The most frequent periods of severe pollution occur from 20:00 to 8:00. Seasonal changes and the northeast monsoon, with wind speeds reaching up to 12 $m/s$, significantly impact NMF3. Therefore, we classify NMF3 as transboundary pollution, accounting for 16.58\% of all features.

The monitoring stations with high contributions to pollutant concentrations in NMF4 are HT, CM, and SL, located in Taichung, with elevated pollution levels during periods of heavy traffic flow. Additionally, high pollution concentrations are associated with a wind speed of 4.8 m/s. As a result, NMF4 is classified as domestic pollution, accounting for 16.43\% of all features.
For NMF5, the stations with high pollutant contributions are TF and ML in Miaoli, with high pollution levels also occurring during periods of heavy traffic flow. Due to the significant influence of northeast monsoons and wind speeds reaching up to 7.6 $m/s$, NMF5 is classified as transboundary pollution, accounting for 16.68\% of all features.
In NMF6, the stations with high pollutant contributions are SY and FY, with lower pollution levels observed between 1:00 and 6:00. High pollution is primarily concentrated at a wind speed of 3.9 $m/s$. Consequently, NMF6 is classified as transboundary pollution, accounting for 16.56\% of all features.

In Fig. \ref{hso2}, we observe that the monitoring stations HT and CM make significant contributions to the SO$_{2}$ concentration in NMF1. Fig. \ref{wso2} indicates that the period with the most severe pollution for NMF1 is from 8:00 to 16:00, with high pollution levels concentrated at a wind speed of 4.8 $m/s$, as shown in Fig. \ref{RSO2}. Thus, we classify NMF1 as domestic pollution, accounting for 21.37\% of all features.

For NMF2, the stations with high pollutant contributions are NT, CS, and PL in Nantou, with the most severe pollution occurring from 10:00 to 20:00. This feature is mainly influenced by westerly winds, with high pollution concentrated at a wind speed of 4 $m/s$. Therefore, NMF2 is classified as domestic pollution, accounting for 24.46\% of all features.

In NMF3, the stations with high contributions to SO$_{2}$ concentrations are HH and EL in Changhua, with the most severe pollution occurring between 8:00 AM and 4:00 PM. Seasonal variations and northeast monsoons, with wind speeds reaching up to 12 $m/s$, significantly impact NMF3. As a result, we classify NMF3 as transboundary pollution, accounting for 17.6\% of all features.

The monitoring station with the highest contribution to SO$_{2}$ concentration in NMF4 is CH. NMF4 is strongly affected by seasonal changes, with higher pollution levels in the spring and September. Additionally, north winds and northeast monsoons influence this feature significantly, with wind speeds of up to 8.1 $m/s$. Therefore, NMF4 is classified as transboundary pollution, accounting for 15.97\% of all features. 

NMF5 has high contributions from the TF station and differs from other features as it shows higher pollution levels during the summer. It is notably affected by strong southwest winds, with wind speeds of up to 7.6 $m/s$. Consequently, we classify NMF5 as transboundary pollution, accounting for 20.61\% of all features.

\subsection*{Analysis on PM10}

In Fig. \ref{hpm10}, we observe that the monitoring stations that contribute significantly to the PM$_{10}$ concentrations in NMF1 are NT, CS, and PL, all located in Nantou. Fig. \ref{wpm10} indicates that the period with the highest pollution levels for NMF1 is from 4:00 PM to 2:00 AM, with concentrations peaking at a wind speed of 4 m/s, as shown in Fig. \ref{RPM10}. Thus, we classify NMF1 as domestic pollution, accounting for 25.22\% of all features.
For NMF2, the highest pollutant contributions come from stations HH, CH, and EL in Changhua. Seasonal changes and northeast monsoons notably impact NMF2, with elevated winter and March pollution levels. Wind speeds can reach up to 9.8 $m/s$, bringing high concentrations of contaminants with the wind. Therefore, we classify NMF2 as transboundary pollution, accounting for 26.18\% of all features. 

The monitoring stations with significant contributions to PM$_{10}$ concentrations in NMF3 are TL, CM, and FY in Taichung. The period with the highest pollution levels for NMF3 is from 10:00 to 1:00. This feature is less affected by seasonal changes, with peak pollution occurring at a wind speed of 3.9 $m/s$. As a result, we classify NMF3 as domestic pollution, accounting for 25.73\% of all features. In NMF4, the stations contributing most to the pollutant concentration are TL, ML, and SY, located in Miaoli. Seasonal variations clearly influence this feature, with higher pollution levels in March and April. Additionally, east winds significantly affect NMF4, with wind speeds reaching up to $7.6 m/s$. Therefore, we classify NMF4 as transboundary pollution, accounting for 22.87\% of all features.

\subsection*{Analysis on PM2.5}

In Fig. \ref{hpm25}, we observe that the monitoring stations with lower contributions to PM2.5 concentrations in NMF1 are TL, CM, FY, NT, CS, and PL. Fig. \ref{wpm10} reveals distinct seasonal variations, with peak pollution levels occurring between 8:00 and 14:00. Additionally, pollution is strongly influenced by the northeast monsoon, as shown in Fig. \ref{RPM25}. Thus, we classify NMF1 as transboundary pollution, accounting for 35.35\% of all features.

High contributions to pollutant concentrations in NMF2 are found at NT, CS, and PL, all situated in the Nantou area. The period of most severe pollution is from 16:00 to 2:00, with March being the peak month. The highest concentrations occur when wind speeds are below 4 m/s. As a result, we classify NMF2 as domestic pollution, accounting for 29.94\% of all features. 
NMF3 contributes more to PM2.5 concentrations at stations SY, HT, TL, CM, and FY. Pollution levels tend to be lower between 4:00 and 8:00, while March and April experience more severe pollution. High concentrations occur when wind speeds are below 5.1 $m/s$. Therefore, we classify NMF3 as domestic pollution, accounting for 34.71\% of all features.  

\section{Method validation}
Chen et al. \cite{D} provide data on the proportion of air pollution in Taiwan for the year 2010. They utilized data from coastal air monitoring stations along with specific meteorological conditions, including concentrations of transboundary pollution. They calculated the total concentration of pollutants by applying the Thiessen polygon area-weighting method. Consequently, the concentration of domestic pollution was derived by subtracting the concentration of transboundary pollution from the total concentration. Their findings indicate that the domestic pollution ratios for SO$_{2}$, NO$_{2}$, and O$_{3}$ are 27\%, 70\%, and 25\%, respectively. We will use these results to validate the findings of our NMF analysis and will explain the experimental results and methodology in the following sections.

\subsection*{Analysis on NO$_{2}$}
\begin{figure}[ph]
	\begin{center}
		\begin{center}
			\includegraphics[height=7cm]{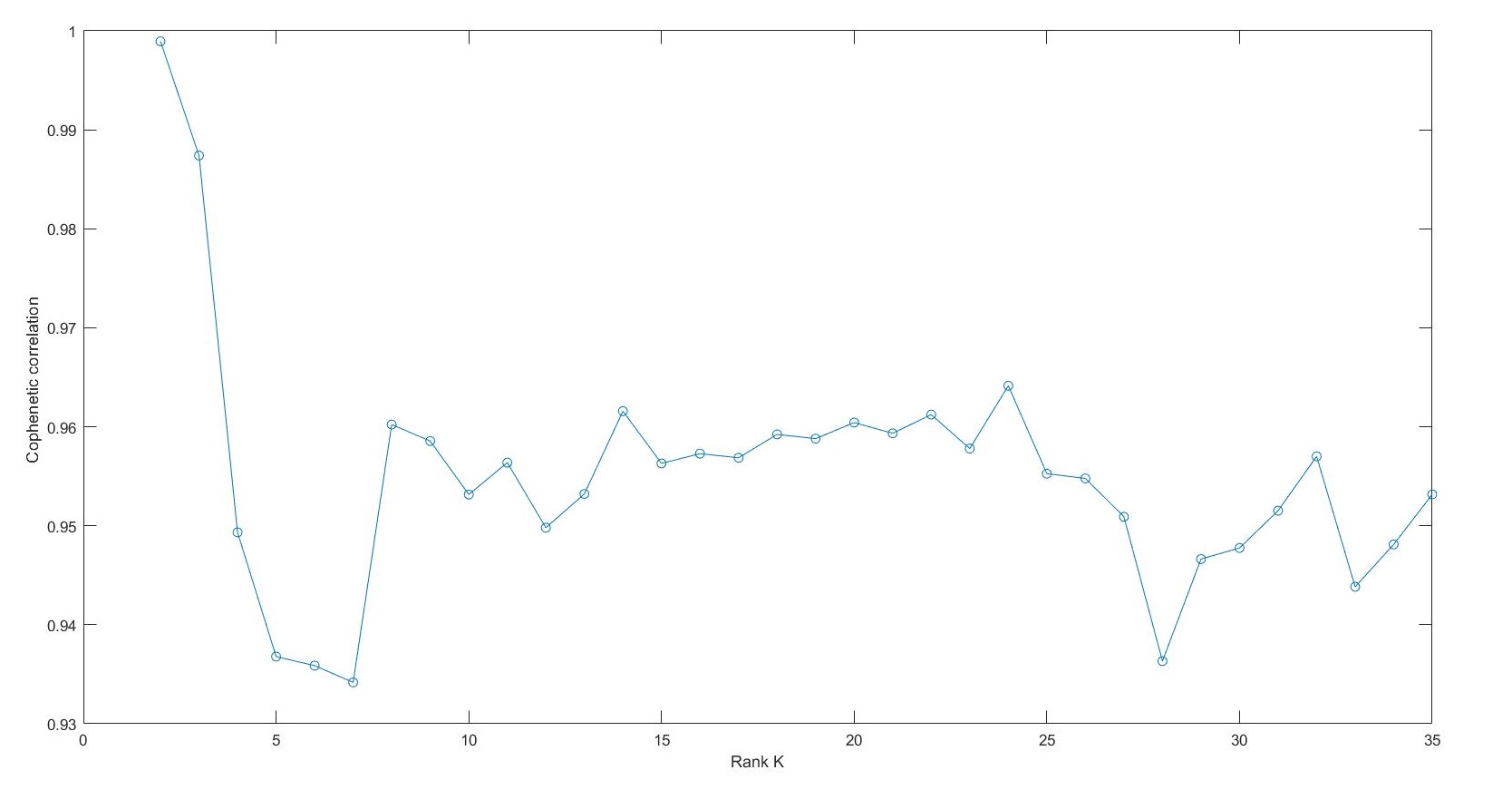}
		\end{center}
		\caption{The percentage of NO$_{2}$: eight components of H}
		\label{Ckno2}
	\end{center}
\end{figure}

\begin{figure}[ph]
	\begin{center}
		\begin{center}
			\includegraphics[height=7cm]{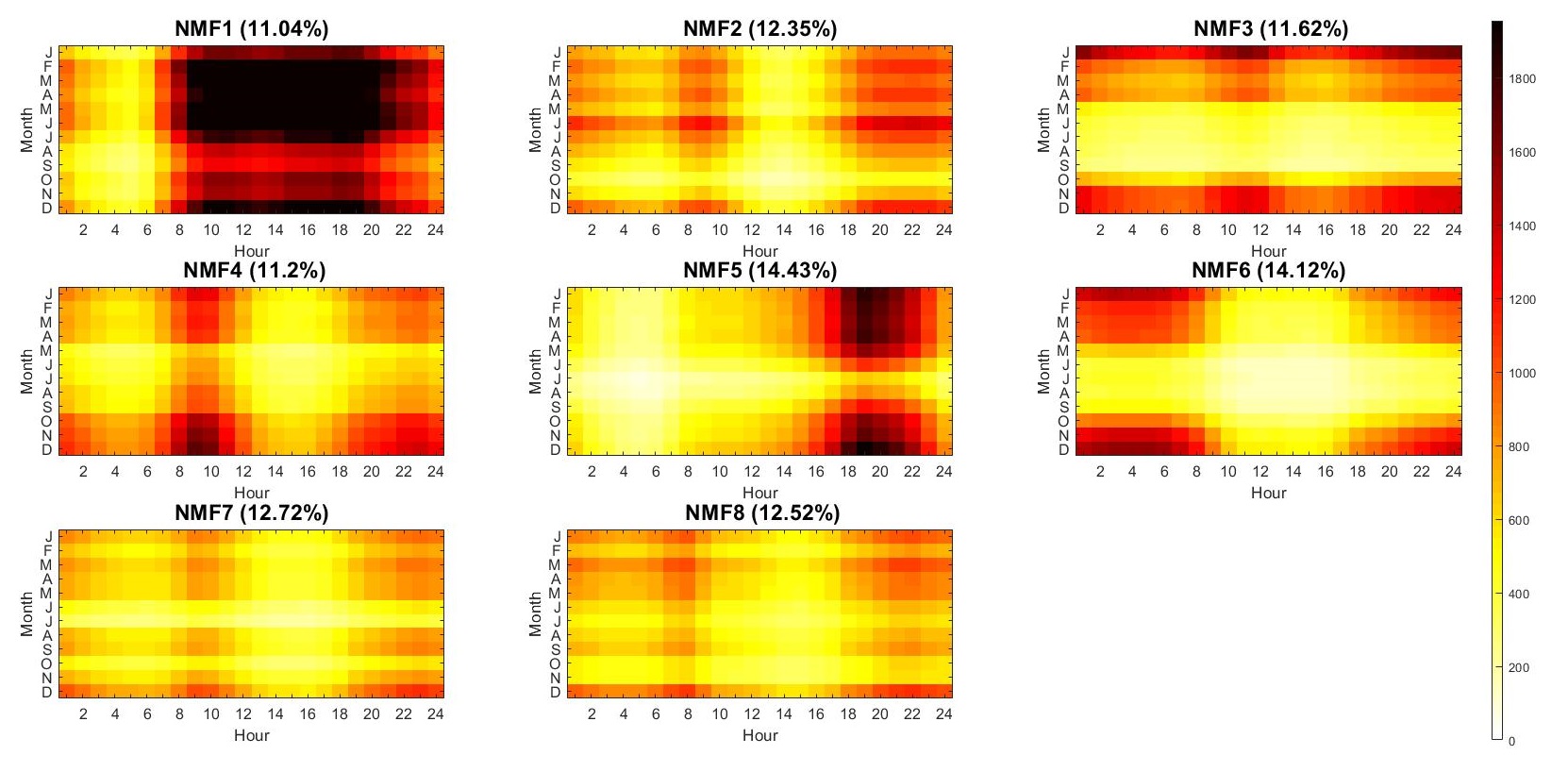}
		\end{center}
		\caption{Day-to-year variation of NO$_{2}$: eight components of W}
		\label{Cwno2}
	\end{center}
\end{figure}

\begin{figure}[ph]
	\begin{center}
		\begin{center}
			\includegraphics[height=12cm]{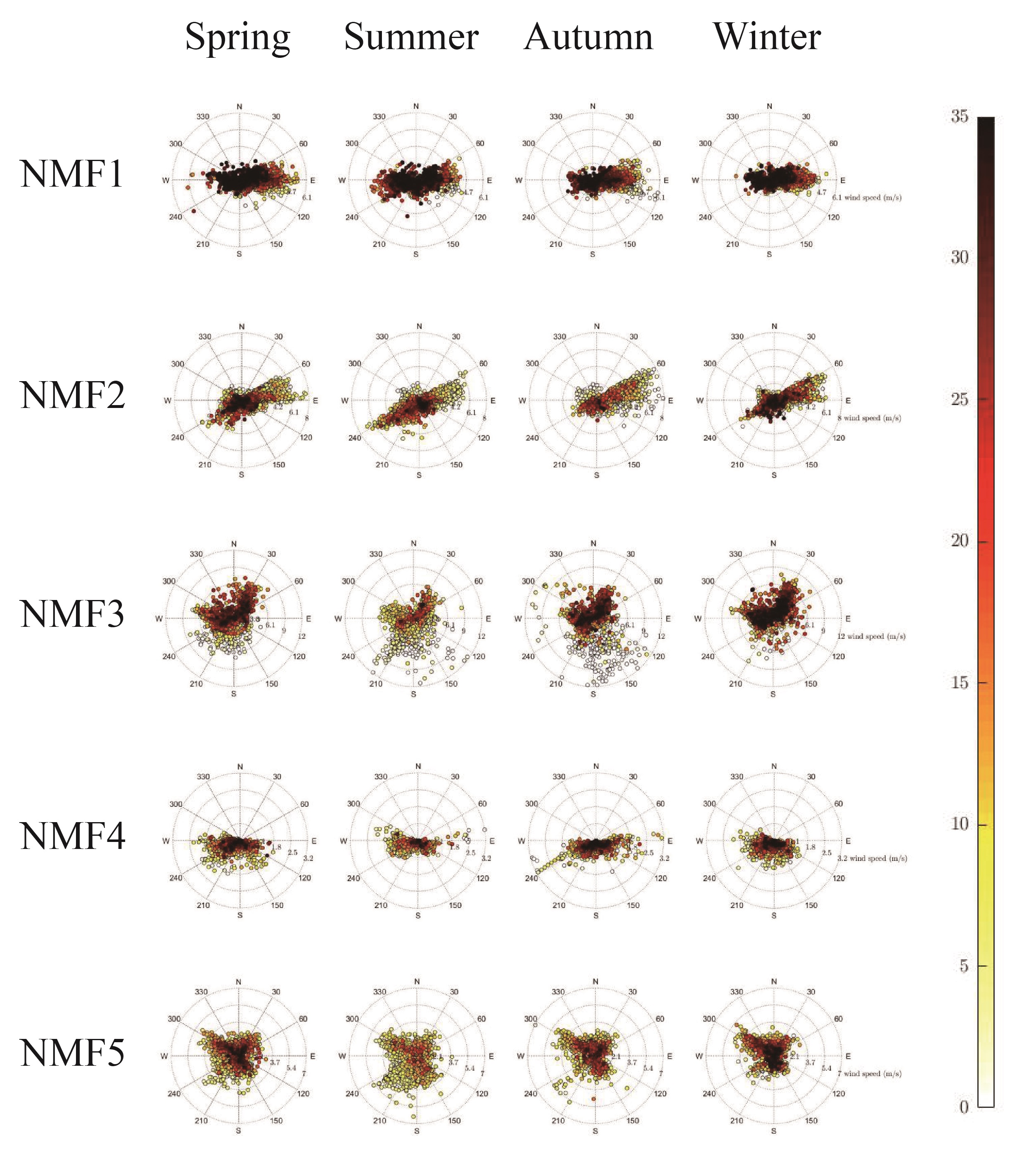}
		\end{center}
		\caption{NMF1-NMF5 for NO$_{2}$ pollution Windrose (the wind direction, wind speed and pollution concentration)}
		\label{Cpno2}
	\end{center}
\end{figure}

\begin{figure}[ph]
	\begin{center}
		\begin{center}
			\includegraphics[height=6cm]{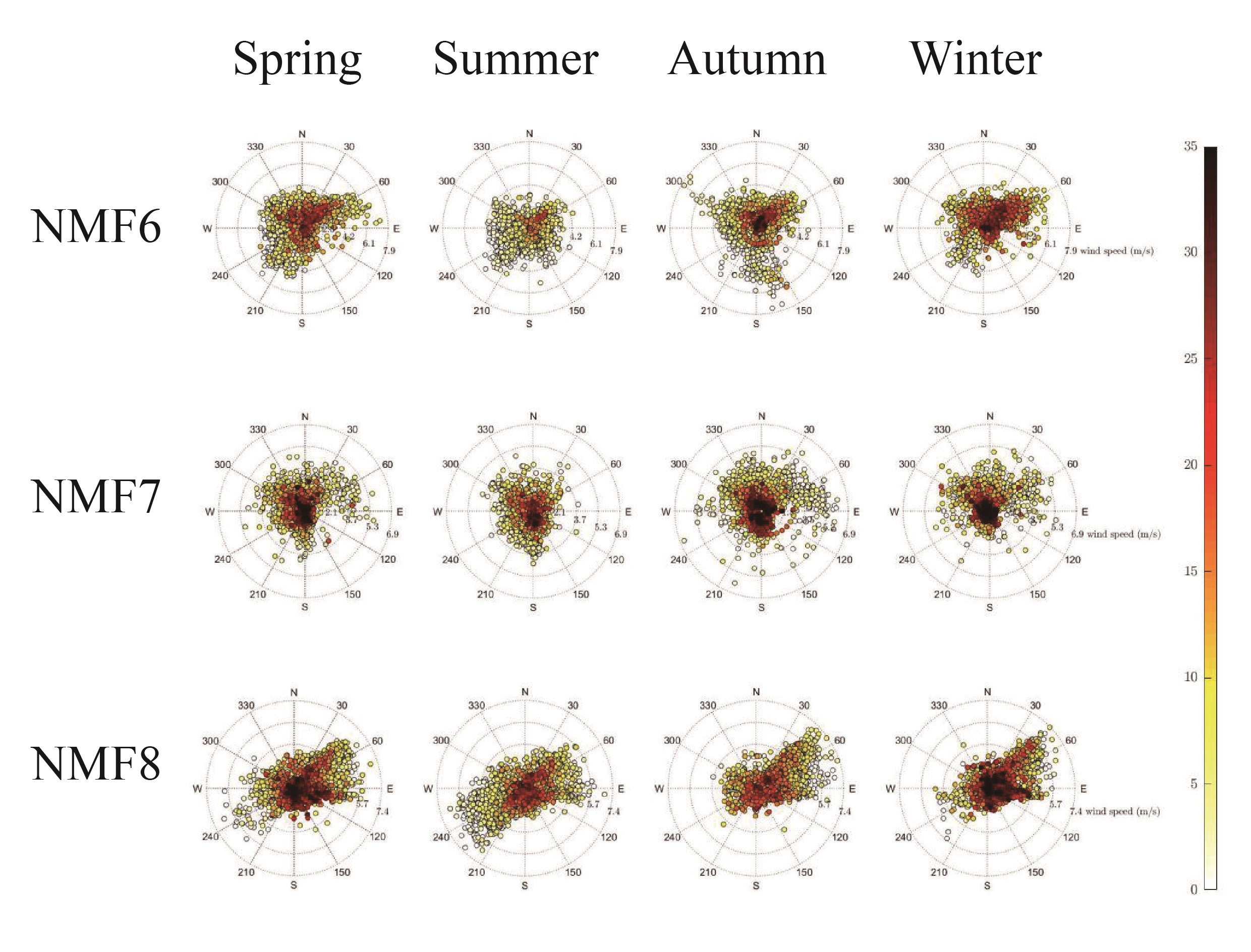}
		\end{center}
		\caption{NMF6-NMF8 for NO$_{2}$ pollution Windrose (the wind direction, wind speed and pollution concentration)}
		\label{Cpno22}
	\end{center}
\end{figure}

We commence our analysis by employing Non-negative Matrix Factorization (NMF) to examine the nitrogen dioxide (NO$_{2}$) pollution data and to determine an appropriate value for the parameter $k$. In this context, we select $k = 8$, as illustrated in Fig. \ref{Ckno2}. Moreover, Figs. \ref{Cwno2}, \ref{Cpno2}, and \ref{Cpno22} demonstrate that the first factor of NMF, labeled NMF1, represents a significant source of NO$_{2}$ pollution, predominantly occurring at wind speeds of 4.7 meters per second. Consequently, we categorize NMF1 as indicative of domestic pollution. 

For NMF2, although wind speeds reach up to 6.1 $m/s$, the highest pollution concentrations occur at wind speeds of 4.2 $m/s$ during times of heavy traffic, leading us to determine that NMF2 is also domestic pollution. NMF3 experiences pollution at wind speeds up to 6.1 $m/s$ and is significantly influenced by seasonal changes. Consequently, we categorize NMF3 as transboundary pollution.
NMF4 to NMF7's primary pollution sources are linked to wind speeds of 2.5 $m/s$, 5.4 $m/s$, 4.2 $m/s$, and 3.7 $m/s$, indicating domestic pollution. NMF8 at 7.4 $m/s$ suggests transboundary pollution. In 2010, transboundary NO$_{2}$ pollution was 24.14\%, while domestic pollution was 75.86\%.

\subsection*{Analysis on O$_{3}$}
\begin{figure}[ph]
	\begin{center}
		\begin{center}
			\includegraphics[height=7cm]{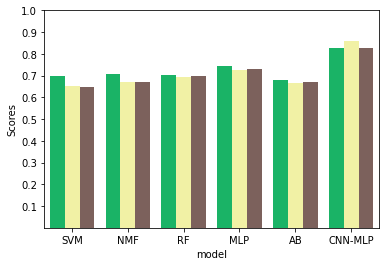}
		\end{center}
		\caption{The percentage of O$_{3}$: ten components of H}
		\label{Cko3}
	\end{center}
\end{figure}

\begin{figure}[ph]
	\begin{center}
		\begin{center}
			\includegraphics[height=7cm]{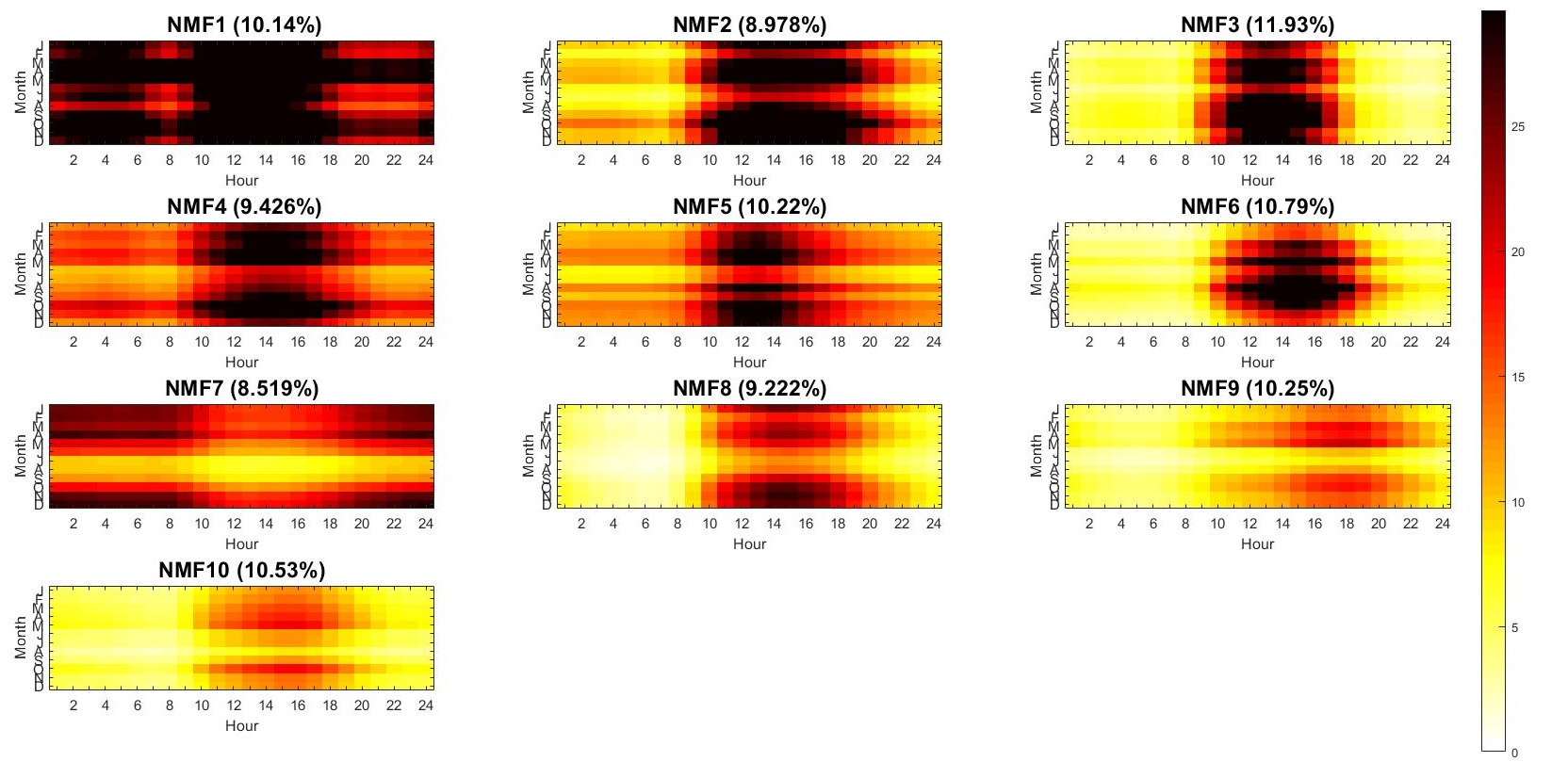}
		\end{center}
		\caption{Day-to-year variation of O$_{3}$: ten components of W}
		\label{Cwo3}
	\end{center}
\end{figure}

\begin{figure}[ph]
	\begin{center}
		\begin{center}
			\includegraphics[height=12cm]{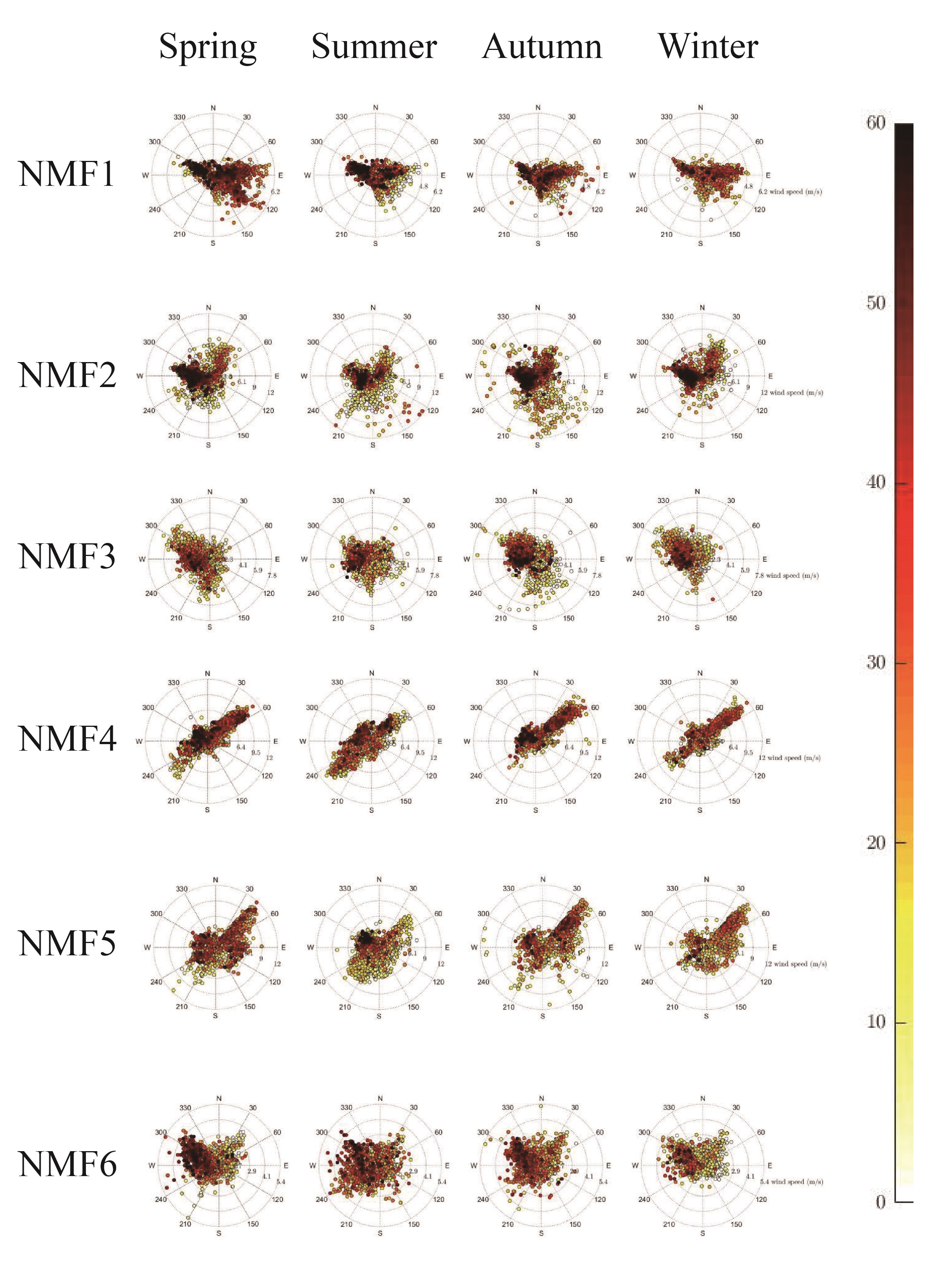}
		\end{center}
		\caption{NMF1-NMF6 for O$_{3}$ pollution Windrose (the wind direction, wind speed and pollution concentration)}
		\label{Cpo3}
	\end{center}
\end{figure}

\begin{figure}[ph]
	\begin{center}
		\begin{center}
			\includegraphics[height=8cm]{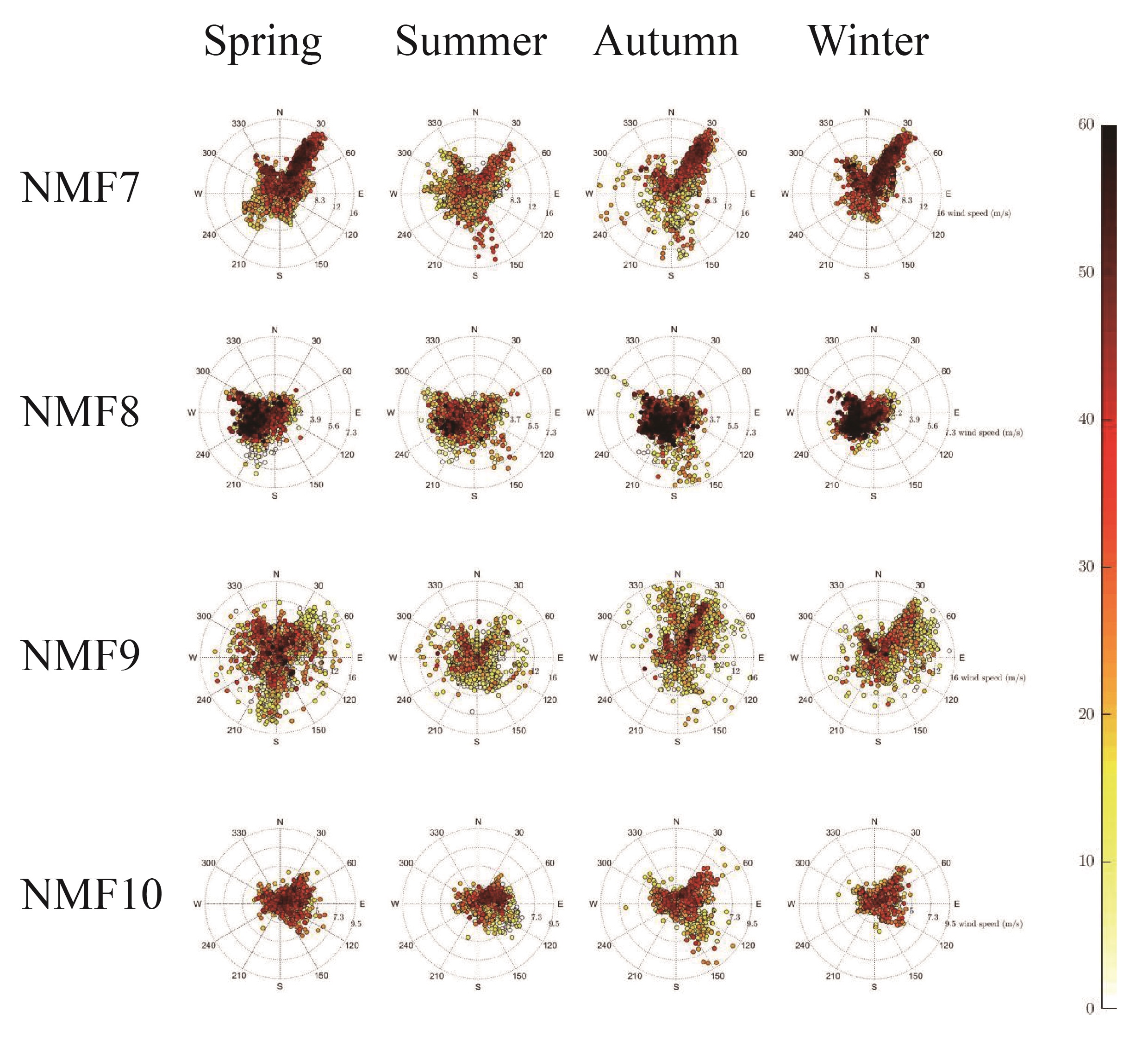}
		\end{center}
		\caption{NMF7-NMF10 for O$_{3}$ pollution Windrose (the wind direction, wind speed and pollution concentration)}
		\label{Cpo32}
	\end{center}
\end{figure}

Subsequently, we conduct an analysis of O$_{3}$ pollution, as depicted in Fig. \ref{Cko3}, wherein we select a value of \( k = 10 \). Figs. \ref{Cwo3}, \ref{Cpo3}, and \ref{Cpo32} demonstrate that the wind speeds for NMF1 and NMF2 are observed at 6.2 $m/s$ and 6.1 $m/s$, respectively. Consequently, we categorize both NMF1 and NMF2 as instances of transboundary pollution. 
NMF3 exhibits high levels of pollution, concentrated at a wind speed of 4.1 $m/s$, significantly influenced by temporal variations, leading to its classification as domestic pollution. Both NMF4 and NMF5 present wind speeds of 12 $m/s$, which compels us to conclude that they, too, constitute transboundary pollution. 

NMF6 reveals significant pollution at a wind speed of 4.1 $m/s$ and is similarly affected by temporal changes; thus, we classify it as domestic pollution as well. Lastly, NMF7 through NMF10 demonstrate wind speeds of 16 $m/s$, 5.6 $m/s$, 16 $m/s$, and 7.3 $m/s$, respectively, all indicating transboundary pollution. In conclusion, for the year 2010, domestic O$_{3}$ pollution accounts for 22.72\%, while transboundary pollution constitutes 77.28\%.

\subsection*{Analysis on SO$_{2}$}
\begin{figure}[ph]
	\begin{center}
		\begin{center}
			\includegraphics[height=7cm]{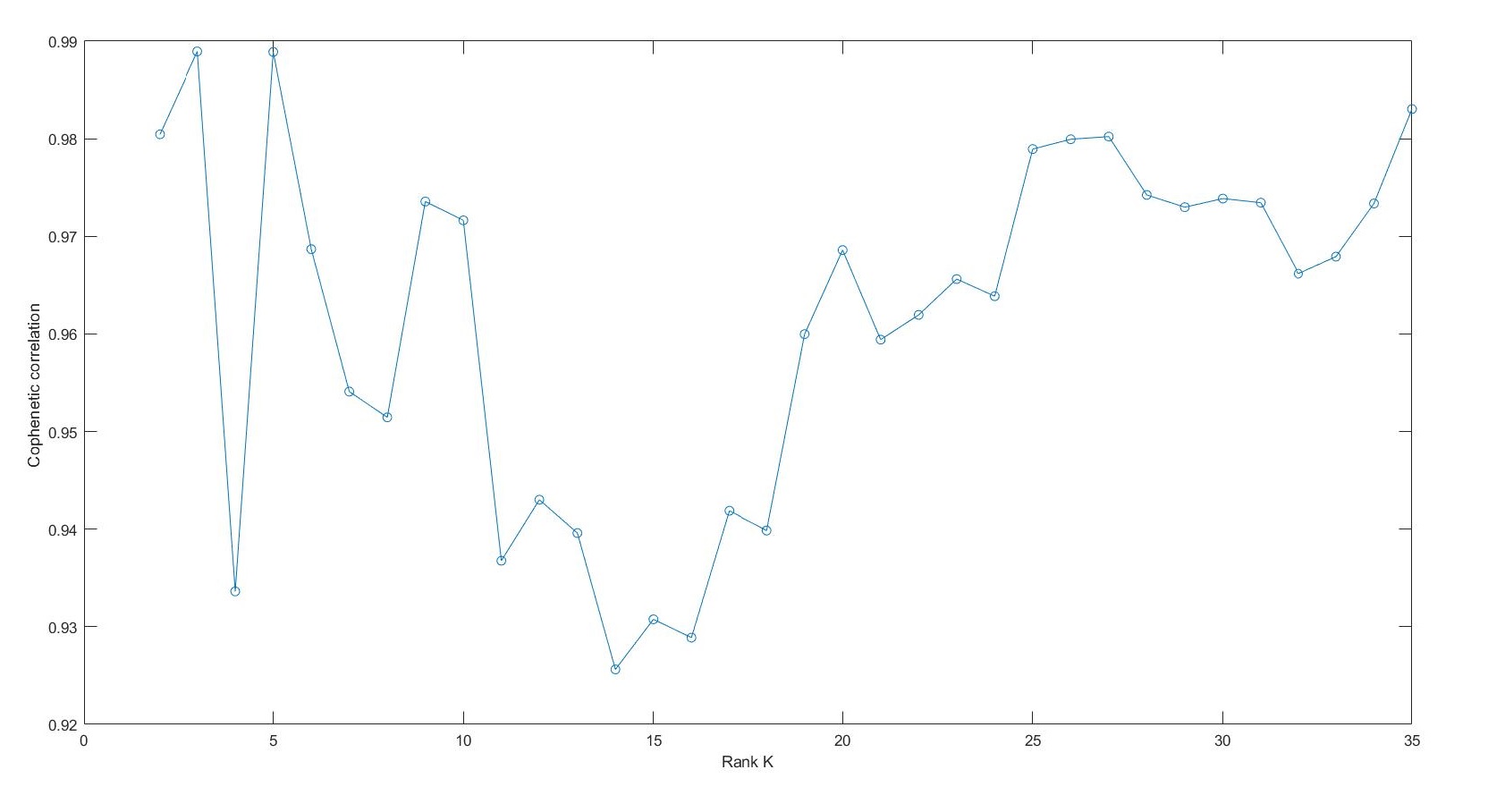}
		\end{center}
		\caption{The percentage of SO$_{2}$: ten components of H}
		\label{Ckso2}
	\end{center}
\end{figure}

\begin{figure}[ph]
	\begin{center}
		\begin{center}
			\includegraphics[height=7cm]{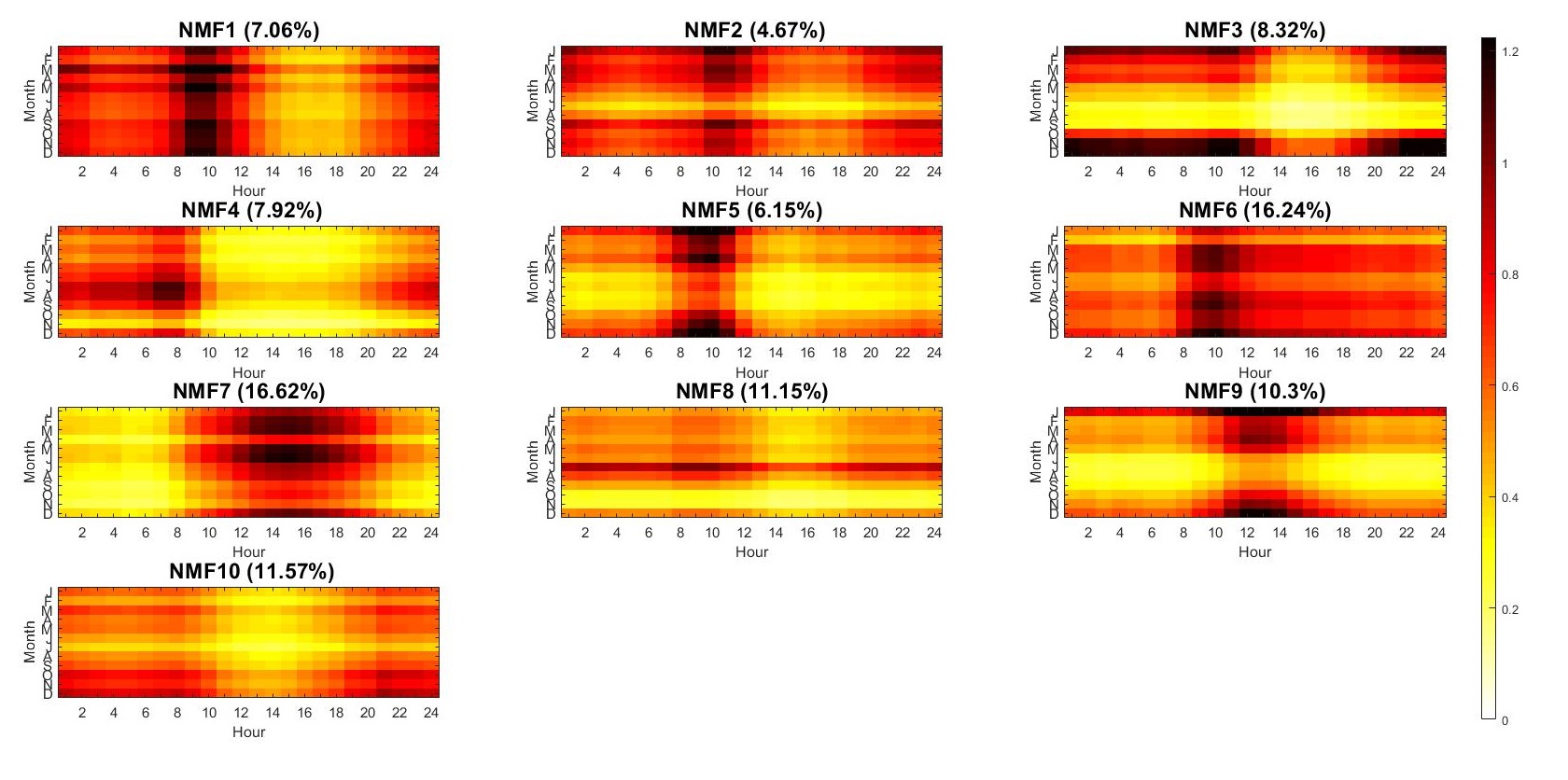}
		\end{center}
		\caption{Day-to-year variation of SO$_{2}$: ten components of W}
		\label{Cwso2}
	\end{center}
\end{figure}

\begin{figure}[ph]
	\begin{center}
		\begin{center}
			\includegraphics[height=10cm]{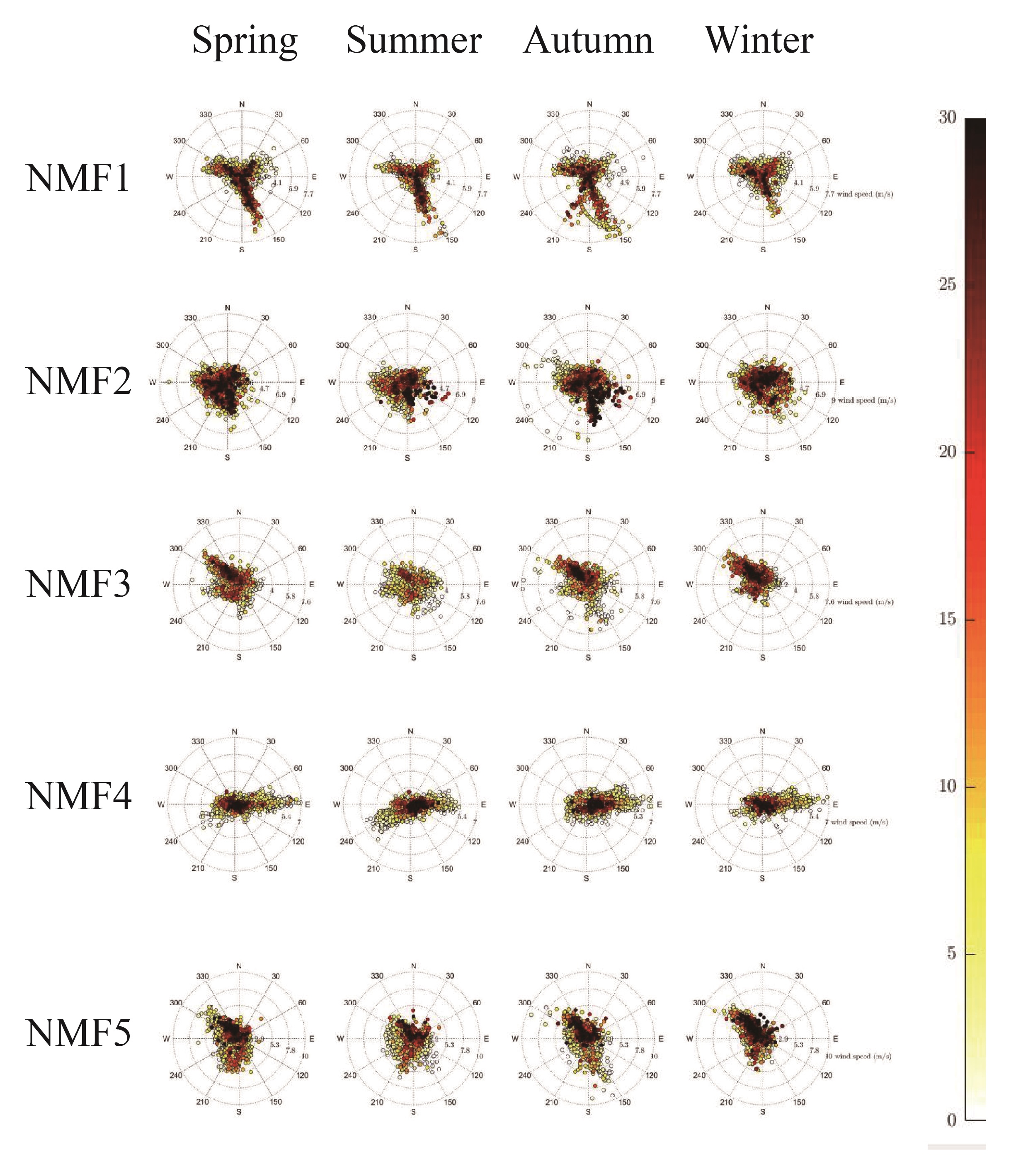}
		\end{center}
		\caption{NMF1-NMF5 for SO$_{2}$ pollution Windrose (the wind direction, wind speed and pollution concentration)}
		\label{Cso2}
	\end{center}
\end{figure}

\begin{figure}[ph]
	\begin{center}
		\begin{center}
			\includegraphics[height=10cm]{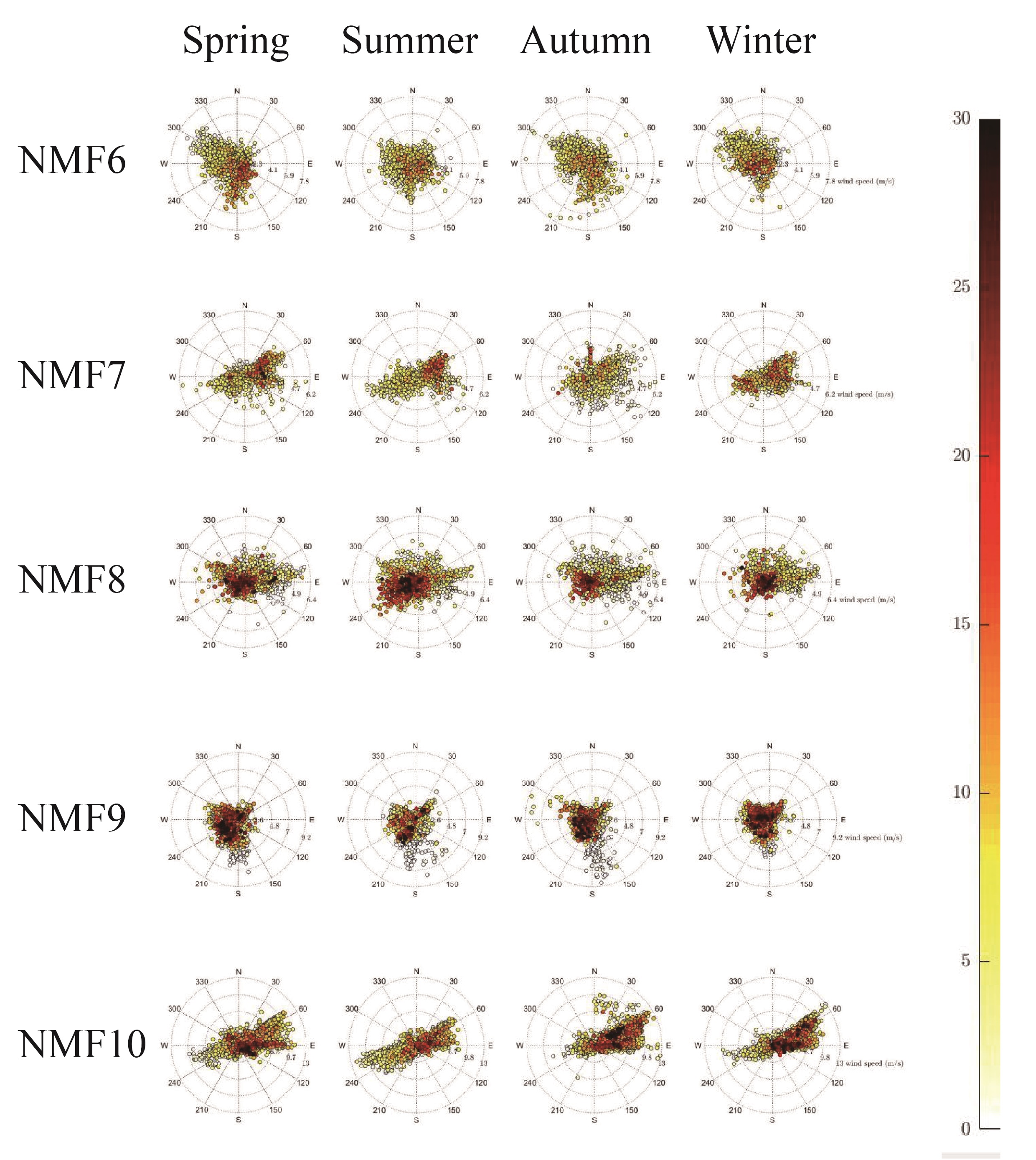}
		\end{center}
		\caption{NMF6-NMF10 for SO$_{2}$ pollution Windrose (the wind direction, wind speed and pollution concentration)}
		\label{Cso22}
	\end{center}
\end{figure}

In conclusion, we conduct an analysis of SO$_{2}$. We have chosen not to select \( $k = 3$ \) or \( $k = 5$ \) in order to preserve the integrity of the data, opting instead for \( $k = 10$ \), as illustrated in Fig. \ref{Ckso2}. From the observations presented in Figs. \ref{Cwso2}, \ref{Cso2}, and \ref{Cso22}, it is evident that NMF1 demonstrates wind speeds reaching up to 7.7 $m/s$, which necessitates its classification as transboundary pollution. 

For NMF2, the distinction between domestic and transboundary pollution remains indeterminate, as it is significantly influenced by seasonal variations, with pollution peaks observed during the autumn season. Consequently, we categorize NMF2 as transboundary pollution. NMF3 to NMF6 exhibit wind speeds of up to 5.8 $m/s$, 7 $m/s$, 7.8 $m/s$, and 5.9 $m/s$, respectively, qualifying them as transboundary pollution. Conversely, NMF7, with a wind speed of $4.7 m/s$, leads us to conclude that it constitutes domestic pollution. NMF8, characterized by wind speeds reaching 6.4 m/s, is also markedly influenced by seasonal changes, resulting in its classification as transboundary pollution. NMF9 reveals elevated pollution concentrations at wind speeds of 4.8 $m/s$, prompting us to classify it as domestic pollution. Lastly, NMF10, exhibiting wind speeds of up to 13 $m/s$, is considered transboundary pollution.
In conclusion, for the year 2010, domestic SO$_2$ pollution represented 26.92\%, while transboundary pollution accounted for 73.08\%. 

\section{Conclusion}
In this study, we employed the Non-negative Matrix Factorization (NMF) method to diminish the dimensionality of air pollution data. Subsequently, we conducted an analysis of wind direction and speed to categorize each pollutant feature as either arising from domestic sources or transboundary origins. Table \ref{4} illustrates the percentages of contributions from domestic and transboundary sources for nitrogen dioxide (NO$_{2}$), sulfur dioxide (SO$_{2}$), particulate matter (PM$_{10}$), and particulate matter (PM$_{2.5}$). The findings indicate that NO$_{2}$ and PM$_{2.5}$ pollution is predominantly influenced by sources within the domestic sphere, whereas transboundary sources significantly impact SO$_{2}$ pollution. For PM$_{10}$, the contributions from domestic and transboundary sources are nearly equivalent. Furthermore, we compared our results with those of Chen et al. \cite{D} as presented in Table \ref{5}. Considering that our findings fall within a $\,\pm 6\%$ margin of deviation from those reported by Chen et al., we deduce that our methodologies are both effective and yield promising outcomes.
\begin{table}[H]
	\begin{center}
		\caption{Proportion of domestic and transboundary pollution.}
		\label{4}
		\begin{tabular}{|c|c|c|}
		\hline
proportion & transboundary pollution ratio & domestic pollution ratio  \\ \hline
 NO$_{2}$  & 33.64\% & 66.36\%	\\ \hline
 SO$_{2}$  & 54.17\% &	45.83\%	\\ \hline
 PM$_{10}$  & 49.05\%	& 50.95\%	 \\ \hline
 PM$_{2.5}$  & 35.35\% &	64.65\%	  \\ \hline
		\end{tabular}
	\end{center}
\end{table}

\begin{table}[H]
	\begin{center}
		\caption{Method validation for 2010 data.}
		\label{5}
		\begin{tabular}{|c|c|c|c|}
		\hline
contaminant & Chen  \cite{D} & NMF  \\ \hline
 NO$_{2}$  & 70\% & 75.9\%	\\ \hline
 SO$_{2}$  & 27\% &	26.9\%	\\ \hline
 O$_{3}$  & 25\%	& 22.7\%	 \\ \hline
		\end{tabular}
	\end{center}
\end{table}


\newpage
\bgroup
\begin{figure}[ph]
	\begin{center}
		\begin{center}
			\includegraphics[height=7cm]{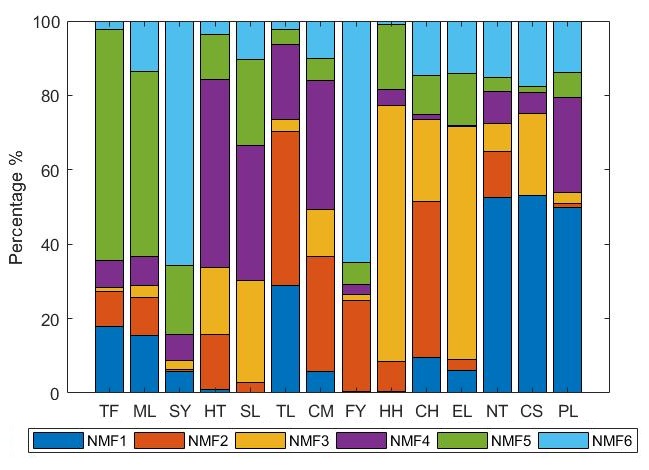}
		\end{center}
		\caption{The percentage of NO$_{2}$: five components of H}
		\label{hno2}
	\end{center}
\end{figure}
\begin{figure}[ph]
	\begin{center}
		\begin{center}
			\includegraphics[height=7cm]{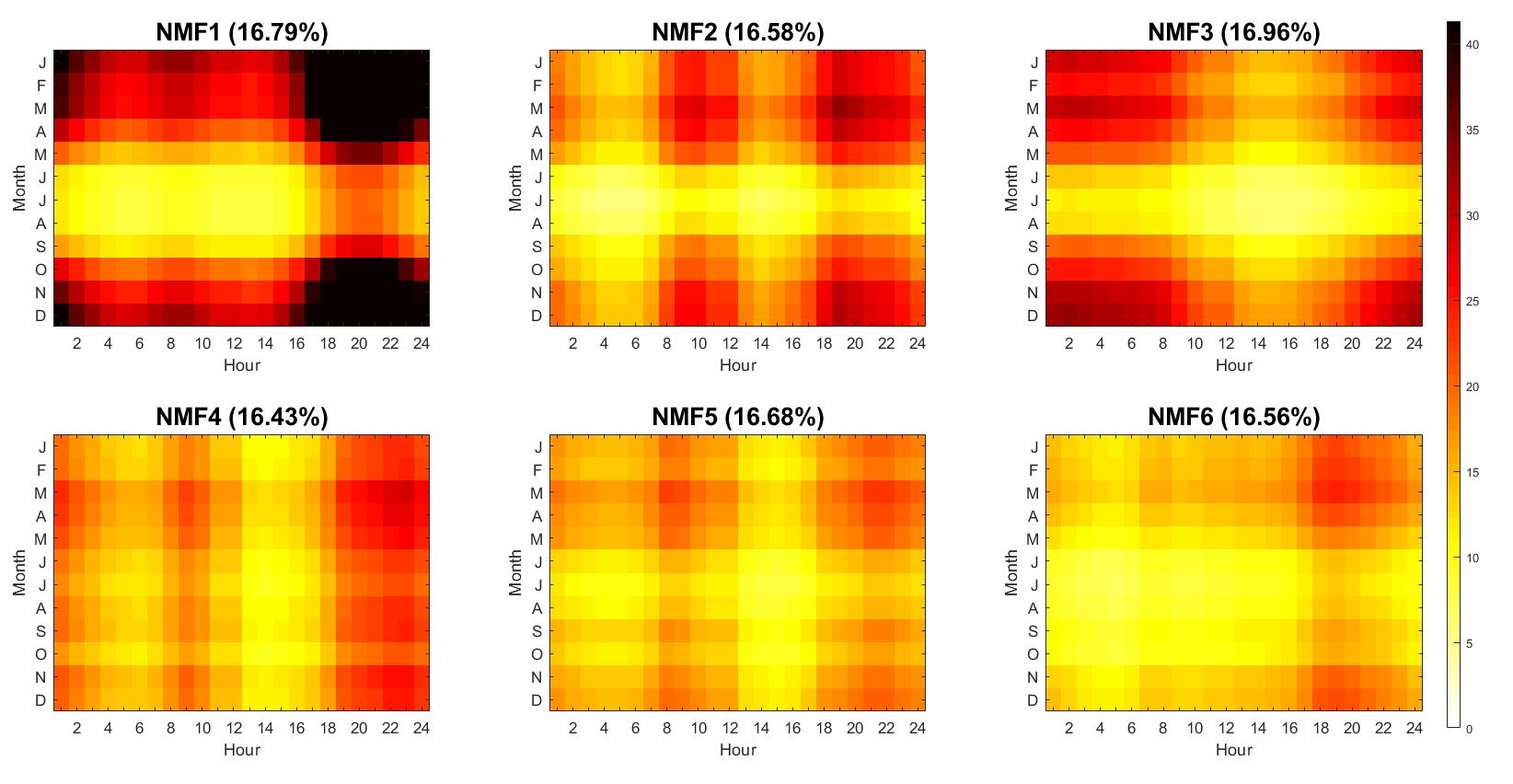}
		\end{center}
		\caption{Day-to-year variation of NO$_{2}$: five components of W}
		\label{wno2}
	\end{center}
\end{figure}

\begin{figure}[ph]
	\begin{center}
		\begin{center}
			\includegraphics[height=12cm]{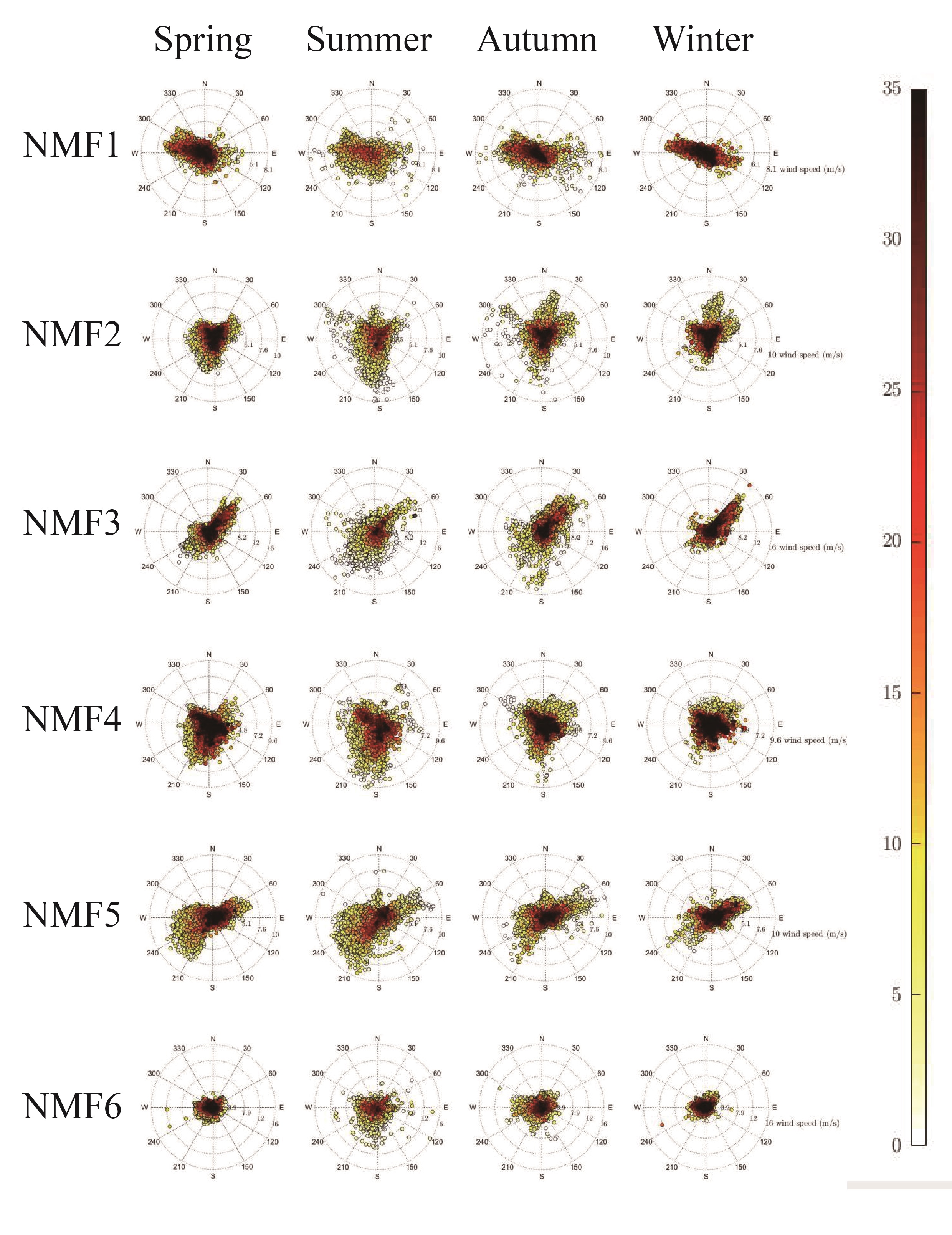}
		\end{center}
		\caption{NMF for NO$_{2}$ pollution Windrose (the wind direction, wind speed, and pollution concentration form a pollution Windrose) }
		\label{RNO2}
	\end{center}
\end{figure}
\egroup

\bgroup
\begin{figure}[ph]
	\begin{center}
		\begin{center}
			\includegraphics[height=7cm]{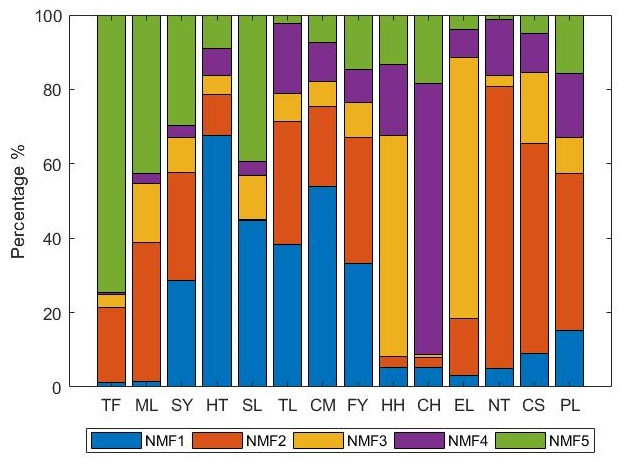}
		\end{center}
		\caption{The percentage of SO$_{2}$: five components of H}
		\label{hso2}
	\end{center}
\end{figure}
\begin{figure}[ph]
	\begin{center}
		\begin{center}
			\includegraphics[height=7cm]{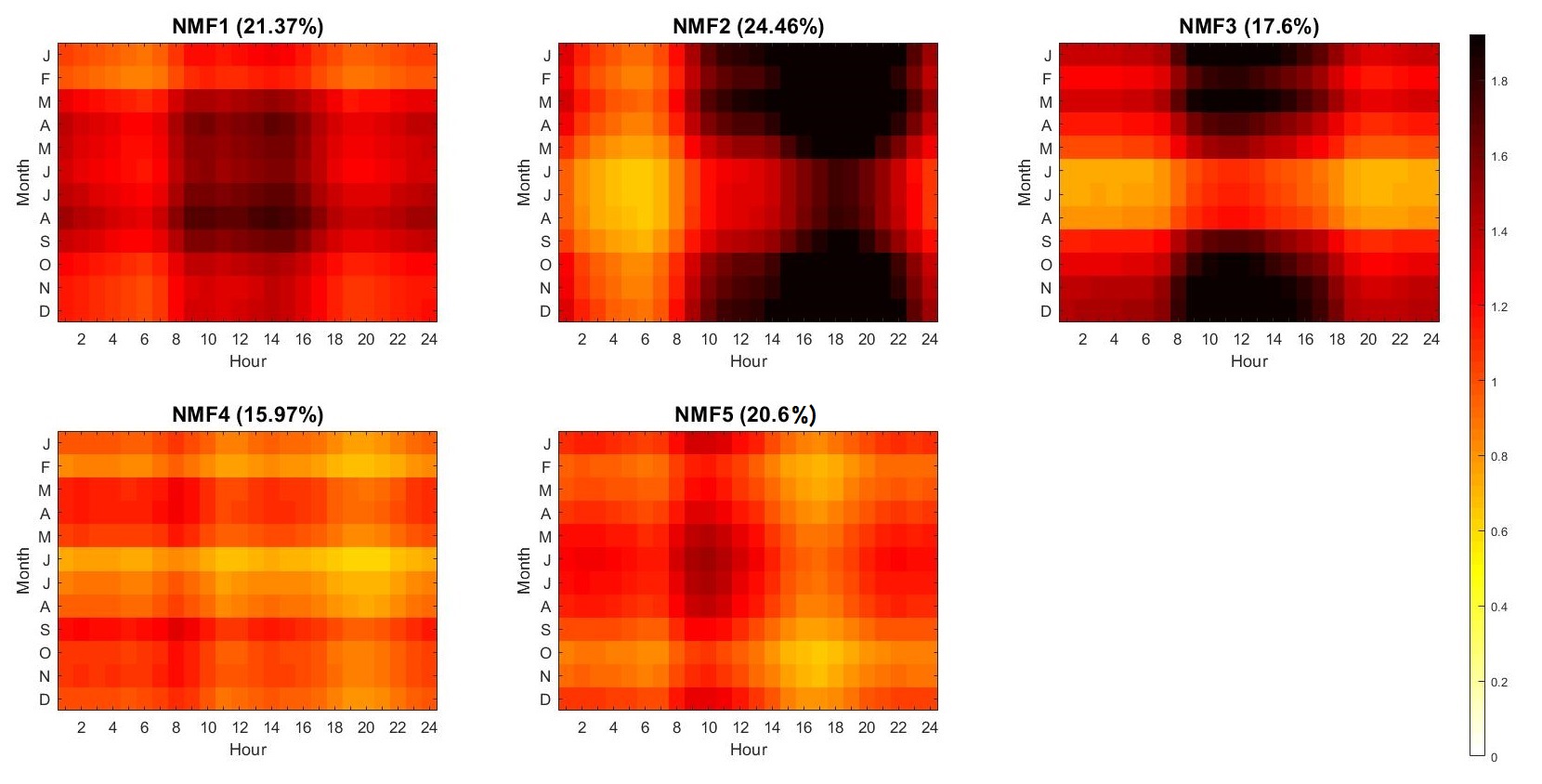}
		\end{center}
		\caption{Day-to-year variation of SO$_{2}$: five components of W}
		\label{wso2}
	\end{center}
\end{figure}

\begin{figure}[ph]
	\begin{center}
		\begin{center}
			\includegraphics[height=10cm]{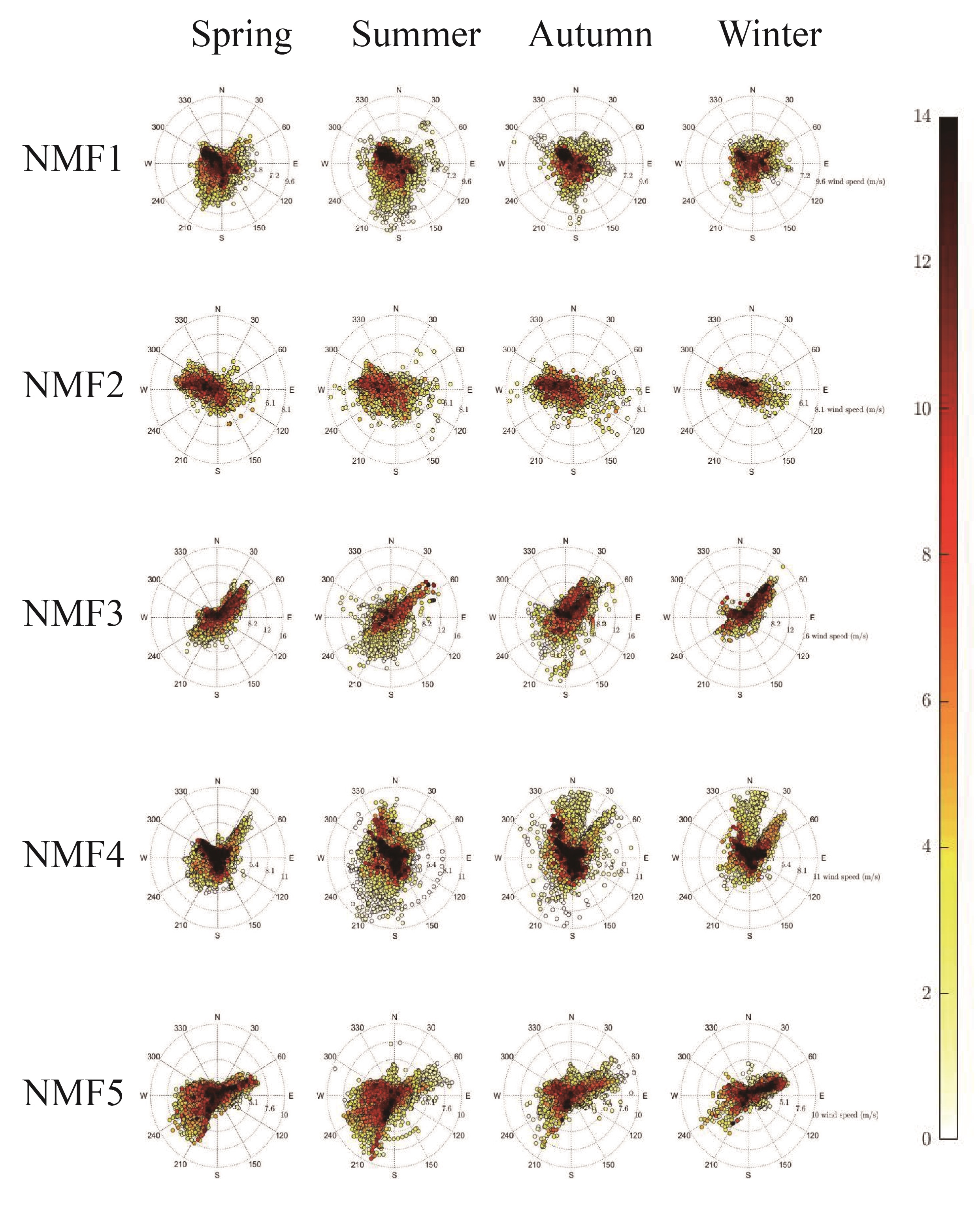}
		\end{center}
		\caption{NMF for SO$_{2}$ pollution Windrose (the wind direction, wind speed, and pollution concentration)}
		\label{RSO2}
	\end{center}
\end{figure}
\egroup

\bgroup
\begin{figure}[ph]
	\begin{center}
		\begin{center}
			\includegraphics[height=7cm]{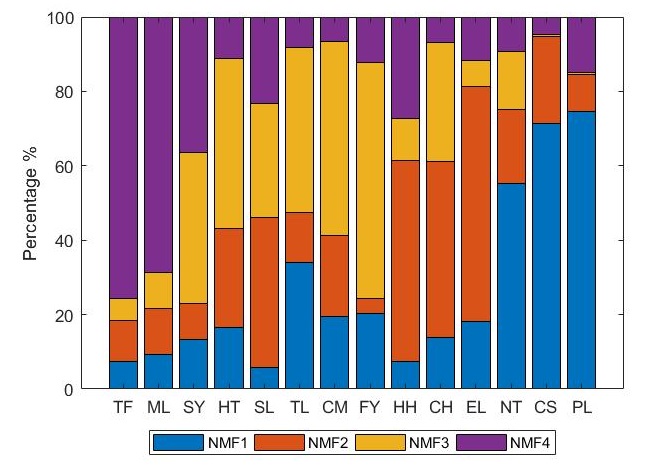}
		\end{center}
		\caption{The percentage of PM$_{10}$: four components of H}
		\label{hpm10}
	\end{center}
\end{figure}
\begin{figure}[ph]
	\begin{center}
		\begin{center}
			\includegraphics[height=7cm]{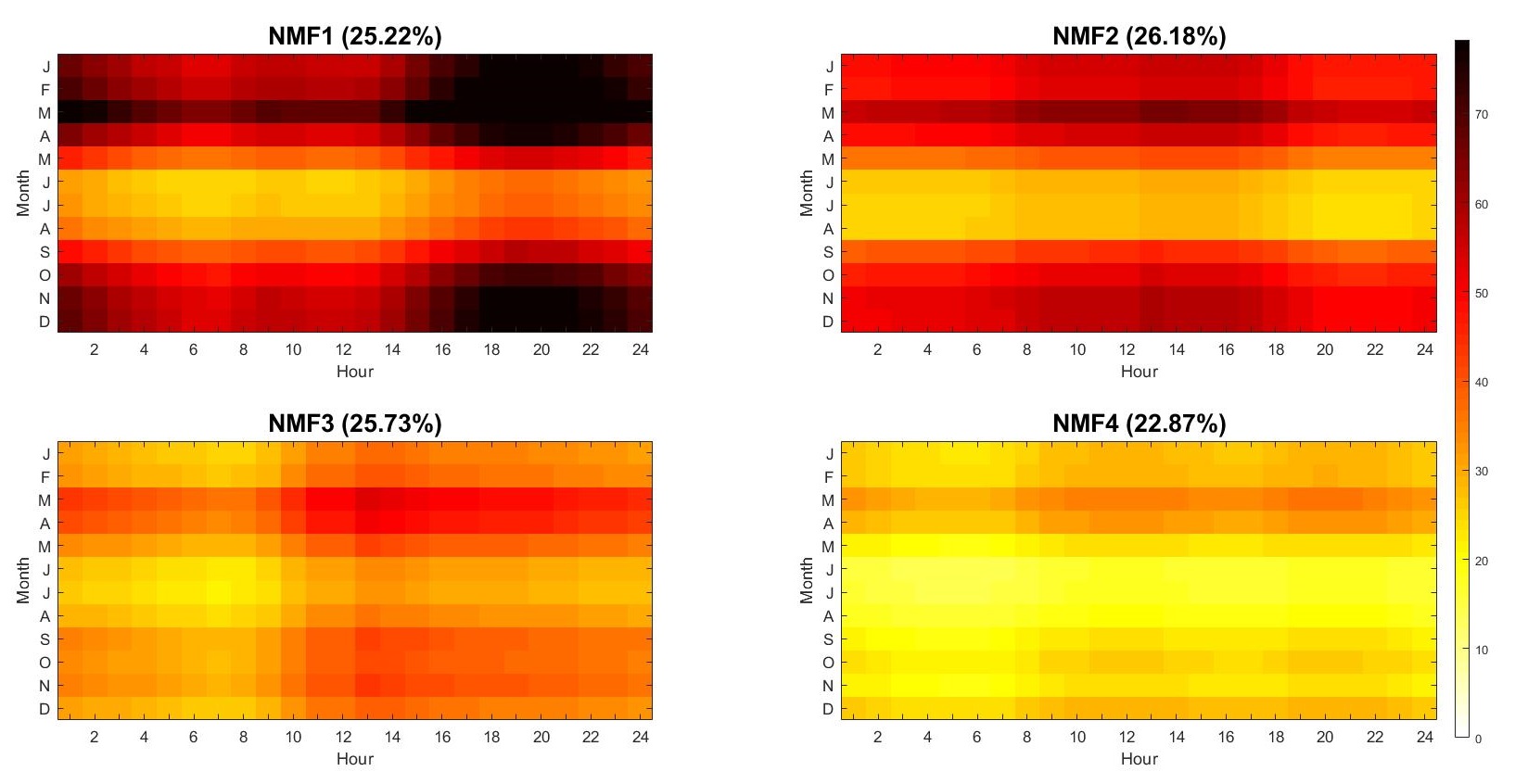}
		\end{center}
		\caption{Day-to-year variation of PM$_{10}$: four components of W}
		\label{wpm10}
	\end{center}
\end{figure}

\begin{figure}[ph]
	\begin{center}
		\begin{center}
			\includegraphics[height=8cm]{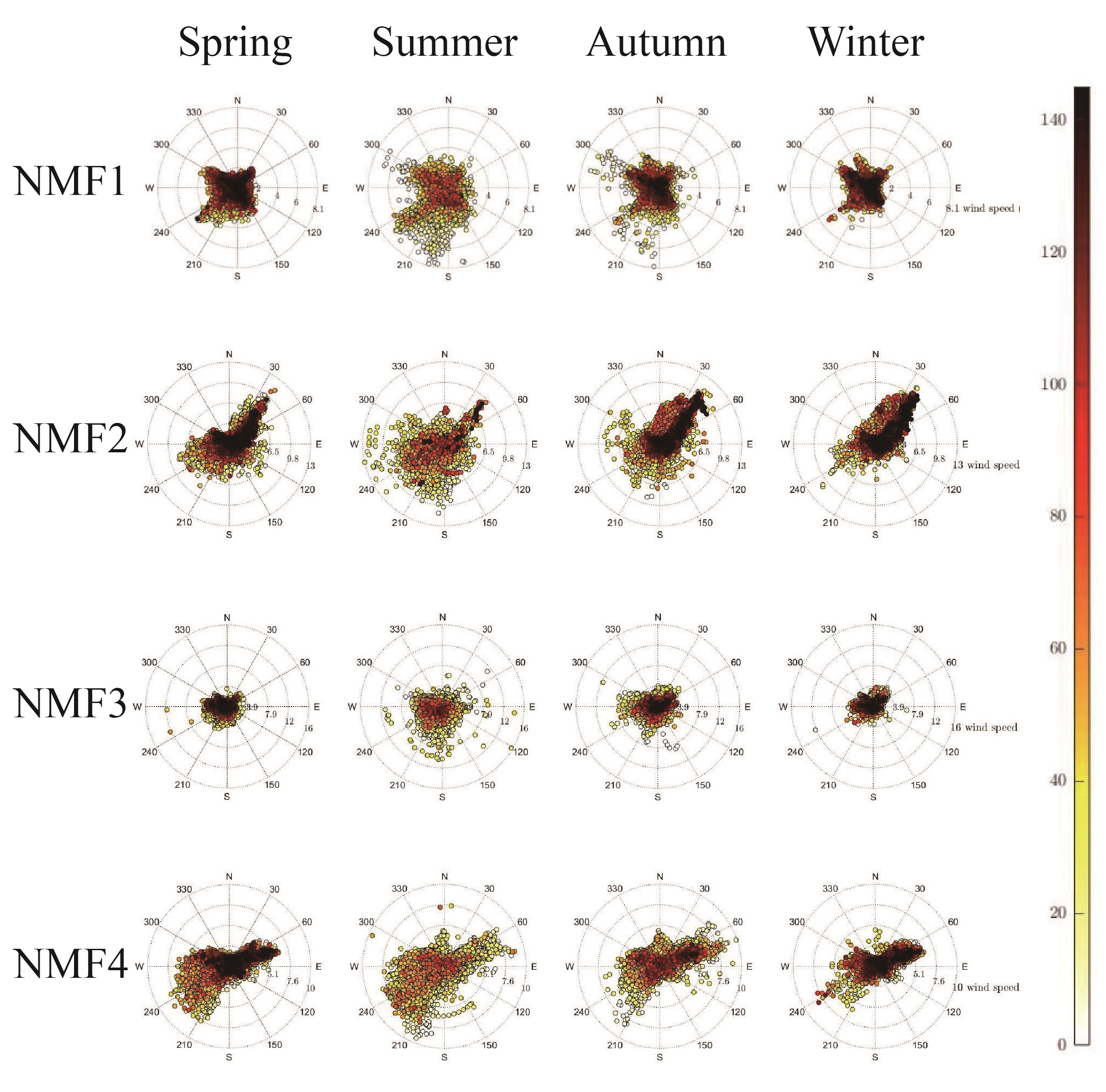}
		\end{center}
		\caption{NMF for PM$_{10}$ pollution Windrose (the wind direction, wind speed, and pollution concentration)}
		\label{RPM10}
	\end{center}
\end{figure}
\egroup

\bgroup
\begin{figure}[ph]
	\begin{center}
		\begin{center}
			\includegraphics[height=7cm]{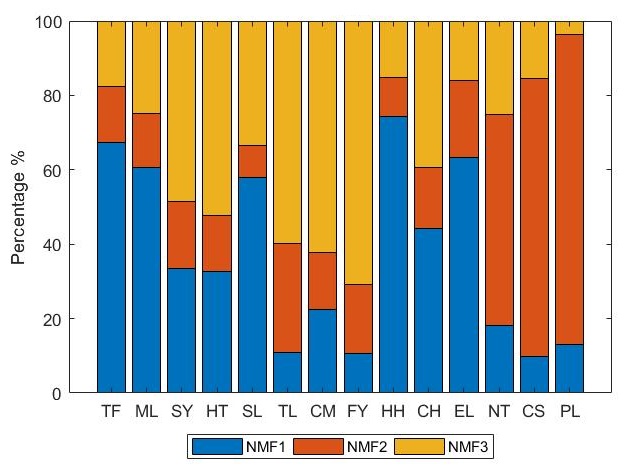}
		\end{center}
		\caption{The percentage of PM$_{2.5}$: three components of H}
		\label{hpm25}
	\end{center}
\end{figure}

\begin{figure}[ph]
	\begin{center}
		\begin{center}
			\includegraphics[height=7cm]{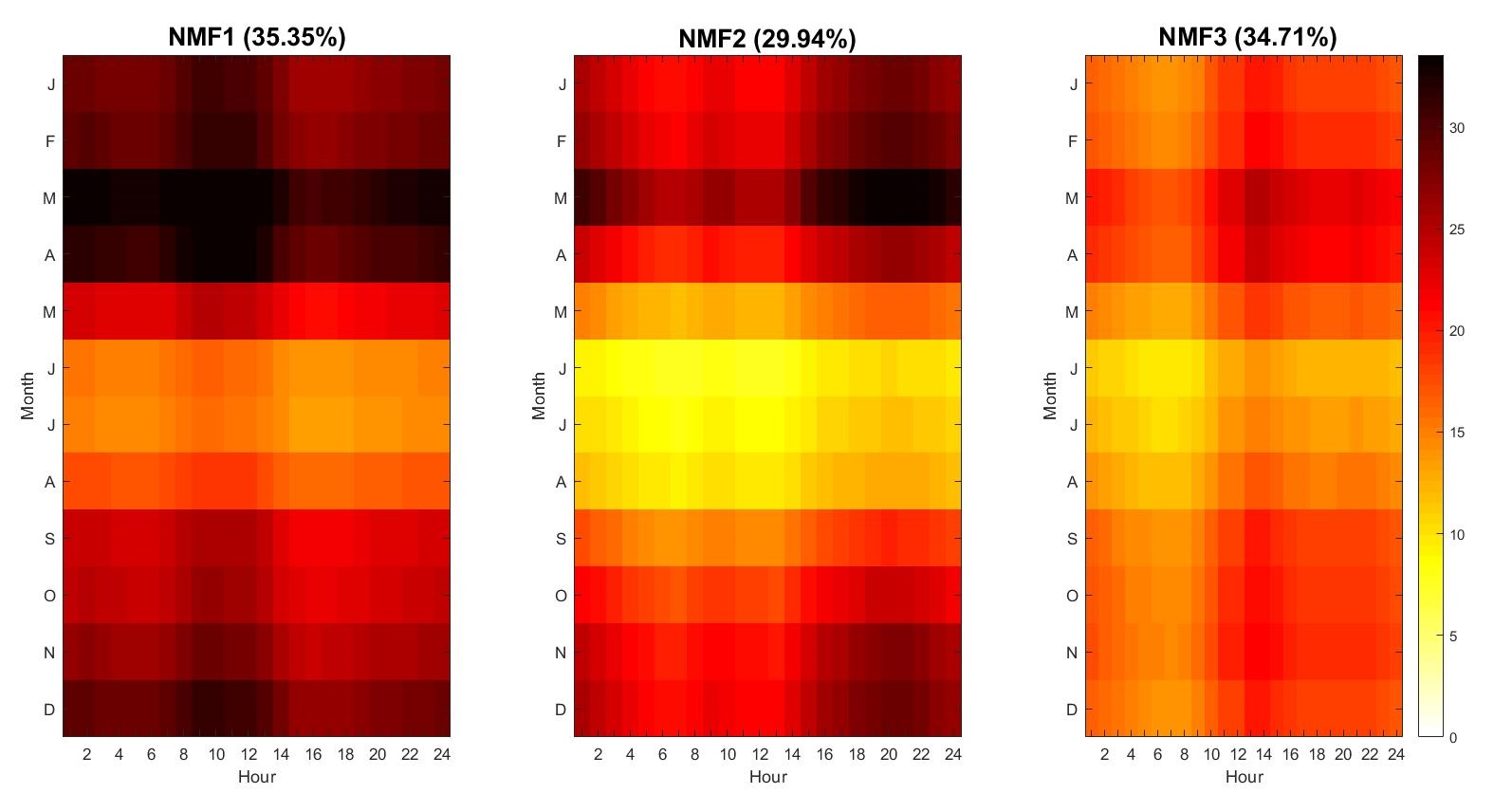}
		\end{center}
		\caption{Day-to-year variation of PM$_{2.5}$: three components of W}
		\label{wpm25}
	\end{center}
\end{figure}

\begin{figure}[ph]
	\begin{center}
		\begin{center}
			\includegraphics[height=7cm]{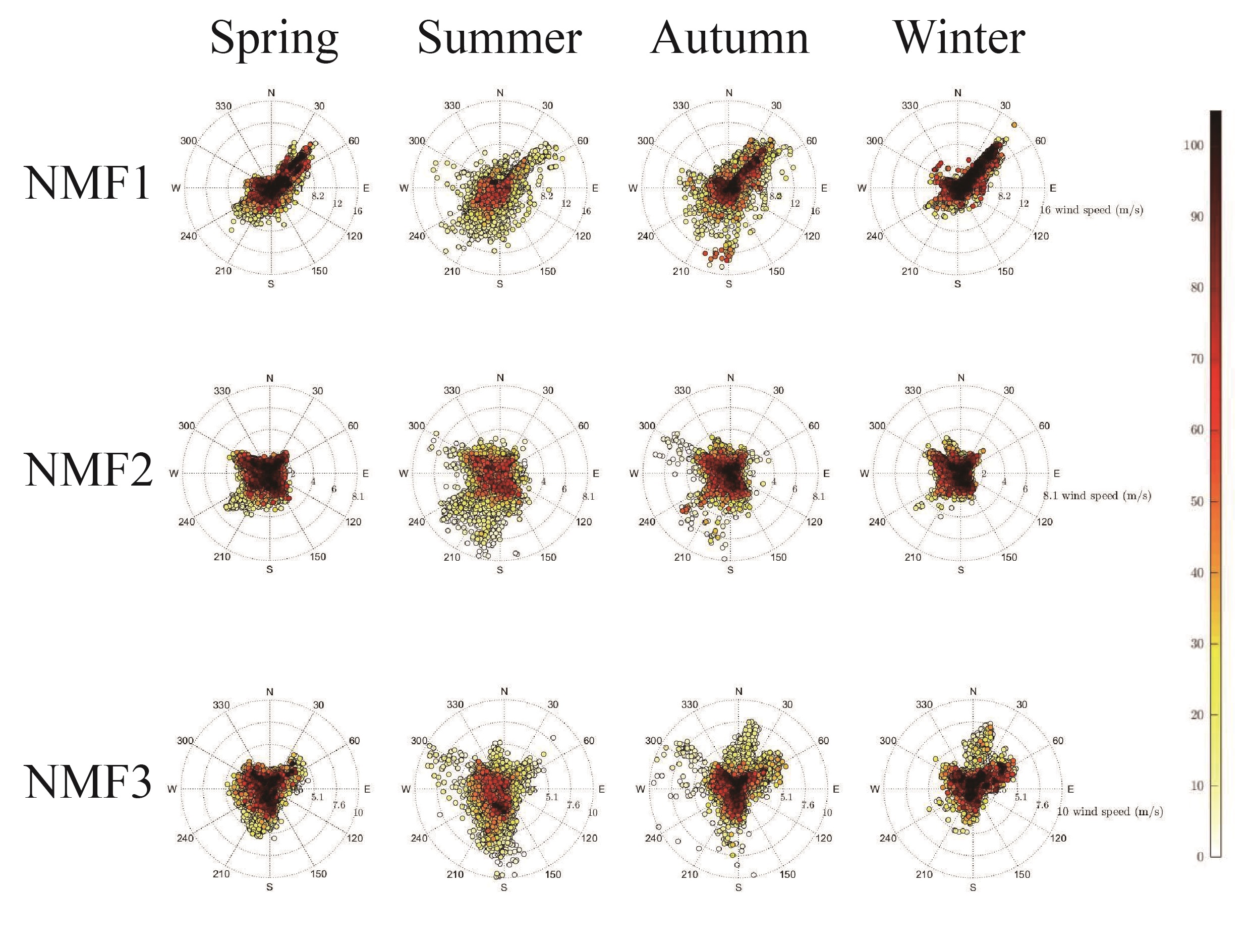}
		\end{center}
		\caption{NMF for PM$_{2.5}$ pollution Windrose (the wind direction, wind speed, and pollution concentration)}
		\label{RPM25}
	\end{center}
\end{figure}
\egroup




	
\end{document}